\documentclass[11pt]{article}
\usepackage[utf8]{inputenc}
\usepackage{amssymb,amsmath}
\usepackage{graphicx,units,mathrsfs,bm}
\usepackage{subfigure}
\usepackage[sort&compress,numbers]{natbib} 
\usepackage[paperwidth=8.5in, paperheight=11in, portrait, top=1in, bottom=1.in, left=1.15in, right=1.15in]{geometry}
\usepackage{tikz}
\usepackage{xcolor}
\usepackage{hyperref}
\usetikzlibrary{calc,decorations.markings,positioning}
\allowdisplaybreaks
\newcommand\bt{\raise 2pt \hbox{$\bigtriangledown$}\hskip 1.5pt}

\newmuskip\pFqskip
\mathchardef\pFcomma=\mathcode`, 

\numberwithin{equation}{section}

\newtheorem{prop}{Proposition}[section]

\newtheorem{defi}{Definition}[section]

\newtheorem{rem}{Remark}[section]

\begin{document}

\title{\bf A unified algebraic underpinning for the Hahn polynomials and rational functions}
\author{
Luc Vinet\textsuperscript{$1$}\footnote{
E-mail: vinet@CRM.UMontreal.CA}~,
Alexei Zhedanov\textsuperscript{$2$}\footnote{
E-mail: zhedanov@yahoo.com} \\[.5em]
\textsuperscript{$1$}\small~Centre de Recherches Math\'ematiques,
Universit\'e de Montr\'eal, \\
\small~P.O. Box 6128, Centre-ville Station, Montr\'eal (Qu\'ebec),
H3C 3J7, Canada.\\[.9em]
\textsuperscript{$2$}\small~School of Mathematics,
Renmin University of China, Beijing, 100872, China
}
\date{\today}
\maketitle

\hrule
\begin{abstract}\noindent

\end{abstract}
An algebra denoted $m\mathfrak{H}$ with three generators is introduced and shown to admit embeddings of the Hahn algebra and the rational Hahn algebra. It has a real version of the deformed Jordan plane as a subalgebra whose connection with Hahn polynomials is established. Representation bases corresponding to eigenvalue or generalized eigenvalue problems involving the generators are considered. Overlaps between these bases are shown to be bispectral orthogonal polynomials or biorthogonal rational functions thereby providing a unified description of these functions based on $m\mathfrak{H}$. Models in terms of differential and difference operators are used to identify explicitly the underlying special functions as Hahn polynomials and rational functions and to determine their characterizations. An embedding of $m\mathfrak{H}$ in $\mathcal{U}(\mathfrak{sl}_2)$ is presented. A Pad\'e approximation table for the binomial function is obtained as a by-product. \\

\bigskip

\hrule


\section{Introduction}
This paper introduces an algebra that subsumes the Hahn algebra \cite{granovskii1992mutual}, \cite{vinet2019heun} and the one recently identified \cite{tsujimoto2020algebraic}, which encodes the bispectral properties of the rational functions of Hahn type. As such it offers a unified algebraic interpretation of both the Hahn polynomials \cite{koekoek2010hypergeometric} and rational functions. This algebra has three generators and derives from a cubic potential \cite{ginzburg2006calabi}. It stands to prove fundamental.

Put simply bispectral problems are situations where functions satisfy a pair of (generalized) eigenvalue equations in the variable and the spectral parameter. Their study was initiated systematically in \cite{duistermaat1986differential} and they have wide-ranging connections and applications \cite{harnad1998bispectral} in particular in time and band limiting \cite{grunbaum2018algebraic}. With their recurrence relations and differential or difference equations, the hypergeometric orthogonal polynomials \cite{koekoek2010hypergeometric} are a particular instance of such problems.

The Askey--Wilson algebra \cite{zhedanov1991hidden} and its relatives such as the Hahn algebra were first obtained as the quadratic algebras realized by the bispectral operators, in the variable or in the degree representations, of the corresponding polynomials. It was observed that their representation theory \cite{zhedanov1991hidden}, \cite{genest2014racah}, \cite{de2015bannai} provides in particular a description of the polynomials whose name they bear. They have moreover been cast in the framework of Leonard pairs \cite{terwilliger2004two} whose classes were shown to correspond to the elements of the Askey scheme.

We shall here pay attention to the Hahn polynomials $Q_n(x, \alpha, \beta, N)$ defined by
\begin{equation}
   Q_n(x, \alpha, \beta, N) = {_3}F_2\left({-n,n+\alpha +\beta +1, -x\atop \alpha +1, -N}; 1 \right ), n=0, 1,\dots N.\label{HahnOP}
\end{equation}
The Hahn algebra $\mathfrak{H}$ which is attached to these polynomials and their duals has two generators $K_1, K_2 $ verifying the relations:
\begin{align}
[K_1,[K_2,K_1]]&=aK_1^2+bK_1+c_1K_2+d_1\mathcal{I} \label{121}\\
[K_2,[K_1,K_2]]&=a\{K_1,K_2\}+bK_2+c_2K_1+d_2\mathcal{I}. \label{212}
\end{align}
As usual $[A,B]=AB-BA$ and $\{A,B\}=AB+BA$. Here $a, b, c_1, c_2, d_1, d_2$ are taken as real parameters. Allowing for affine transformations of the generators their number can be reduced to two. When $\mathfrak{H}$ is realized by the bispectral operators of the polynomials, they are given in terms of the parameters $\alpha$ and $\beta$ of the polynomials. It is known (e.g. from the literature already cited) that the overlaps between eigenvectors of $K_1$ and $K_2$ in
representations of dimension $N+1$ are given in terms of Hahn polynomials.

Biorthogonal rational functions (BRFs) are known to arise as solutions of Generalized Eigenvalue Problems (GEVPs) involving two tridiagonal matrices \cite{zhedanov1999biorthogonal}. It is also understood that some of these BRFs \cite{rahman1981families}, \cite{wilson1991orthogonal}, \cite{ismail1995generalized}, \cite{gupta1998contiguous}, \cite{rahman1993classical}, \cite{spiridonov2003theory} are bispectral in that they are solutions of two GEVPs. A first foray into the algebraic interpretation of the bispectrality of BRFs was recently achieved in \cite{tsujimoto2020algebraic} by looking at the rational functions of Hahn type:
\begin{equation}
   \mathcal{U}_n(x;\alpha,\beta,N)= \frac{(-1)^n(-N)_n}{(\beta+1)_n} \: {_3}F_2\left({-x,-n,\beta+n-N\atop -N,\alpha-x}; 1 \right ), n=0, 1,\dots N. \label{ratHahn}
\end{equation}
(We use the standard notation for Pochhammer symbols and generalized hypergeometric series \cite{koekoek2010hypergeometric}.) This led to the extension of Leonard pairs to Leonard pencils corresponding to the two GEVPs defining the bispectral framework. These two pencils were found to be composed out of three operators $X, Y, Z$ which in the case of the rational functions $\mathcal{U}_n$ given above were found to satisfy the rational Hahn algebra $r\mathfrak{H}$ with the following defining relations:
\begin{align}
    [Z,X] &= Z^2 +Z \label{CMZX} \\
    [X,Y] &= \xi_1(X^2+Z^2) + \{X,Z\} +\{Y,Z\} + \xi_2 X + \xi_3 Z +Y + \xi_0 \mathcal{I} \label{CMXY} \\
    [Y,Z] &= 3 X^2 + Z^2 +\xi_1 \{X,Z\} +\xi_4 X + \xi_2 Z + \xi_0 \mathcal{I}. \label{CMYZ}
\end{align}
The parameters $\xi_i, i=0, 1, 2, 3, 4$ are a priori taken to be arbitrary real constants; these become expressions involving $\alpha$ and $\beta$ in the realization in terms of the operators that define the bispectral problem of which the rational functions $\mathcal{U}_n$ are solutions.

We shall explain below that both the Hahn algebra $\mathfrak{H}$ and the rational Hahn algebra $r \mathfrak{H}$ can be embedded in an algebra with three generators $\{X, Z, V\}$ which we will call the meta-Hahn algebra $m\mathfrak{H}$ and describe in the next section. This is the basis of the joint account of the two bispectral problems.  When dealing with a linear pencil, $P=A+\sigma B$, one can set as was done in \cite{derevyagin2012cmv}, \cite{tsujimoto2017tridiagonal} the two-parameter problem $Pu=\tau u$. \textit{Depending on whether $\sigma$ or $\tau$ is taken as the eigenvalue, one is either looking at a GEVP involving the operators $A-\tau$ and $B$ or at an ordinary eigenvalue problem (EVP) for the operator $P$}. This is how the treatment of the Hahn polynomials and rational functions will be unified and how $m\mathfrak{H}$ will be seen to provide a synthesizing algebraic framework for these two families of special functions.

The exposition will proceed as follows. In addition to presenting $m\mathfrak{H}$ in Section 2 and providing its Casimir element, we shall indicate that it is built as an enlargement of the deformed Jordan plane. How the Hahn algebra $\mathfrak{H}$ and the rational Hahn algebra $r \mathfrak{H}$ inject into the meta-Hahn algebra $m\mathfrak{H}$ will be described in Section 3. Certain elements of the $m\mathfrak{H}$ representation theory will be developed in Section 4. We shall consider bases corresponding to the GEVP defined by the generators $X$ and $Z$ and to the adjoint GEVP, and also bases corresponding to the EVPs associated to $V$ and to the pencil $X+\mu Z$ as well as to their adjoints. The bispectral properties of various overlaps between these bases will be established in Section 5. A differential realization of $m\mathfrak{H}$ will be given in Section 6 and used on the one hand to obtain the explicit expressions of the bispectral overlaps and on the other to derive properties of these special functions. The biorthogonal Hahn rational functions and the Hahn orthogonal polynomials will be recovered there. The embedding of $m\mathfrak{H}$ in $\mathcal{U}(\mathfrak{su}(2))$ is the topic of Section 7. A summary and an outlook will serve as conclusion in Section 8. Two appendices complete the paper. Appendix A will provide computational details on the orthogonality of the eigenvectors of $V$ and of its adjoint $V^T$ in the differential model of $m\frak{H}$. Appendix B looks at Jacobi polynomials with a parameter that is a negative integer. It will offer the restricted versions of Gauss' and Pfaff's transformation formulas that are valid in this case and that were used in Section 6. It will furthermore include a Pad\'e approximation table for the binomial function as an interesting by-product.

\section{The meta-Hahn algebra $m\mathfrak{H}$}

We start with the algebra that has two generators $X$ and $Z$ verifying
\begin{equation}
 [Z,X]=Z^2+Z.  \label{ZXA}
\end{equation}
We will comment below on this noteworthy algebra, see Remark \ref{rem2}. We now extend it by the addition of one extra generator $V$ and a central element $\mathcal{I}$. We assume that the commutation relations of $V$ with $X$ and $Z$ minimally depart from those of Lie algebras by exhibiting a single quadratic term and that they are of the form
\begin{equation}
[X,V] = \eta_6 \{V,Z\}+ \eta_1 X + \eta_7 V + \eta_2 Z + \eta_0 \mathcal{I},\label{XVA}
\end{equation}
\begin{equation}
[V,Z] = \eta_4 X + \eta_5 Z + \eta_3 \mathcal{I},\label{VZA}
\end{equation}
where $\eta_0, \dots \eta_7$ are some (real) parameters. Enforcing the Jacobi relation

\begin{eqnarray}
[V,[Z,X]] + [Z,[X,V]] + [X,[V,Z]]=0 \label{Jac_meta}
\end{eqnarray}
leads to:

\begin{prop}
The relations  \eqref{ZXA} -  \eqref{XVA} are compatible with the Jacobi identity \eqref{Jac_meta} if and only if the constraints
\begin{eqnarray}
\eta_6=\eta_7=1, \quad \eta_5=\eta_1 \label{constr_meta}
\end{eqnarray}
are imposed.

\end{prop}
It is appropriate to fix a standardized form of the emerging algebra so as to identify how many essential parameters it contains. This number is reduced to three with the help of the following automorphisms:
\begin{equation}
    V\rightarrow \kappa V \qquad X\rightarrow X+\zeta Z \label{aut}
\end{equation}
that preserve the commutation relation between $X$ and $Z$.
We shall not consider here the degenerate case $\eta_4 = 0$; under this proviso it is possible to choose $\kappa$ so that $\eta_4=2$ and to pick $\zeta$ in order to have $\eta_2 = -\eta_1$. This brings us to the following:
\begin{defi}
The meta-Hahn algebra $m\mathfrak{H}$ is generated by the three elements $V$, $X$ and $Z$ that satisfy the relations
\begin{equation}
 [Z,X]=Z^2+Z, \label{CZX}
\end{equation}
\begin{equation}
[X,V] = \{V,Z\}+ \eta_1 X + V - \eta_1 Z + \eta_0 \mathcal{I}, \label{CXV}
\end{equation}
\begin{equation}
[V,Z] = 2 X + \eta_1 Z + \eta_3 \mathcal{I}. \label{CVZ}
\end{equation}
\end{defi}

It is further observed that $m\mathfrak{H}$ derives from a potential. Given $F_{\it f}=\mathbb{C}[x_1, x_2, \dots x_n]$ a free associative algebra with $n$ generators, $F_{cycl}=F_{\it f}/[F_{\it f},F_{\it f}]$. $F_{cycl}$ has the cyclic words $[x_{i_1} x_{i_2}\dots x_{i_r}]$ as basis. The cyclic derivative $\frac{\partial}{\partial x_j}: F_{cycl} \rightarrow F_{\it f}$ is such that
\begin{equation}
\frac{\partial [x_{i_1} x_{i_2}\dots x_{i_r}] }{\partial x_j}= \sum_{\{s|i_s=j\}} x_{i_s +1}x_{i_s +2}\dots x_{i_r}x_{i_1}x_{i_2}\dots x_{i_s -1}
\end{equation}
and is extended to  $F_{cycl}$ by linearity. Let $\Phi(x_1,\dots x_n) \in F_{cycl}$. An algebra whose defining relations are given by
\begin{equation}
\frac{\partial\Phi}{\partial x_j}=0, \qquad j=1,\dots n
\end{equation}
is said to derive from the potential $\Phi$. Now let $x_1=V$, $x_2=X$ and $x_3=Z$ and take
\begin{equation}\label{eq:Phi}
\Phi=[XVZ]-[XZV]-[VZ^2]-[VZ]-[X^2]+\frac{\eta_1}{2}[Z^2]-\eta_1[XZ]-\eta_3[X]-\eta_0[Z]
\end{equation}
It is readily seen that the relations \eqref{CZX}, \eqref{CXV}, and \eqref{CVZ} of $m\mathfrak{H}$ are given by
$\frac{\partial \Phi}{\partial V}=0$, $\frac{\partial \Phi}{\partial Z}=0$, $\frac{\partial \Phi}{\partial X}=0$ respectively.

This fact points at the existence of a non-trivial central element with terms corresponding to those in $\Phi$ with the exclusion of $[XVZ]-[XZV]$.
\begin{prop}
The Casimir element of $m\mathfrak{H}$ is given by:
\begin{equation}
Q = \{V, Z^2+Z\} + 2X^2 +  (2-\eta_1) Z^2 +\eta_1 \{X,Z\} + 2 \eta_3  X +2(\eta_0+1) Z. \label{Q}
\end{equation}
\end{prop}
Let us record that before the automorphisms \eqref{aut} are used to set $\eta_2=-\eta_1$ and $\eta_4=2$ the element $Q$ such that $[Q,X]=[Q,Z]=[Q,V]=0$ has the following expression:
\begin{eqnarray}
Q = \{V, Z^2+Z\} + \eta_4 X^2 + \eta_1 \{X,Z\} +(\eta_2 +\eta_4)Z^2 + 2 \eta_3 X +(2\eta_0+ \eta_4)Z  \label{Q_expr}
\end{eqnarray}
which is of course consistent with \eqref{Q}.
\begin{rem}
We shall at times allow for $\eta_0$ and $\eta_1$ to be in fact central.
\end{rem}
\begin{rem}\label{rem2}
The meta-Hahn algebra has a two-dimensional subalgebra generated by $X$ and $Z$ with commutation relation given by \eqref{CZX}. Its remarkable connection to Hahn polynomials will be stressed in Sections 4 and 6. If we consider this algebra over the field of complex numbers, upon setting $X=\frac{i}{2}\hat{X}$ and $Z=\frac{i}{2}\hat{Z}-\frac{1}{2}$, the commutation relation is converted to $[\hat{Z}, \hat{X}]=-\hat{Z}^2 -1$. This is recognized as the deformed Jordan plane denoted $\mathcal{J}_1$ in \cite{gaddis2015two}.
 This algebra is known \cite{gaddis2013pbw} to be a PBW deformation of the two-dimensional Artin-Schelter regular algebra $[X,Z]=Z^2$, which is called the Jordan plane or meromorphic Weyl algebra and has been much studied, see for instance \cite{zhedanov2005regular}, \cite{iyudu2014representation}, \cite{shirikov2018two}. (We would like to also cite \cite{benkart2013parametric} for the study of the representations of the algebras defined by the relation $[X,Z]=h(Z)$ with $h$ a polynomial in $Z$.)
\end{rem}

\section{Embeddings of the Hahn and rational Hahn algebras in the meta-Hahn algebra}

One main feature of $m\mathfrak{H}$ is that it contains the algebras $\mathfrak{H}$ and $r\mathfrak{H}$; this explains the prefix meta.
We now describe  the two embeddings.

\subsection{Hahn algebra: $\mathfrak{H} \hookrightarrow m\mathfrak{H}$}

Introduce the element
\begin{equation}
 W=X + \mu Z \label{W_mu}
\end{equation}
where $\mu$ is an arbitrary real parameter. Using the commutation relations \eqref{CZX}, \eqref{CXV} and \eqref{CVZ} of $m\mathfrak{H}$, it is checked that upon taking
\begin{equation}
    K_1=W \qquad \text{and} \qquad  K_2=V,
\end{equation}
 the relations \eqref{121} and \eqref{212} of the Hahn algebra $\mathfrak{H}$ are satisfied with
\begin{align}
    &a=2, \; b=2\mu-\eta_1 +2\eta_3,\; c_1=-1,\; d_1=-Q,\\
    &c_2=-\eta_1 (\eta_1 +2), \; d_2= \eta_0 (2\mu -\eta_1)-\eta_1 \eta_3 (1+\mu).\label{coeffH}
\end{align}
The following proposition thus follows:
\begin{prop}
The centrally extended Hahn algebra admits an embedding in the meta-Hahn algebra with $V$ and $W$ as generators and with the parameters of \eqref{121} and \eqref{212} given by \eqref{coeffH}.
\end{prop}
The proof is straigthforward. Note that the Casimir element $Q$ of $m\mathfrak{H}$ given in \eqref{Q} occurs in the parameter $d_1$; to be precise, it is hence a central extension of the Hahn algebra that is thus injected in the meta-Hahn algebra.

Observe that the commutation relations between the generators $V, W, Z$ form a closed set that reads:
\begin{equation}
    [Z,W]=Z^2+Z, \label{CrZW}
\end{equation}
\begin{equation}
  [W,V] = \{V,Z\} + \tau_1 W + V + \tau_2 Z + \tau_0 \mathcal{I}, \label{CrWV}
\end{equation}
\begin{equation}
   [V,Z] = 2 W + \tau_1Z +\tau_3\mathcal{I},\label{CrVZ}
\end{equation}
with
\begin{equation}
   \tau_0=\eta_0-\mu \eta_3, \quad \tau_1=\eta_1-2\mu,  \quad \tau_2=2\mu^2-2\eta_1 \mu -\eta_1, \quad \tau_3=\eta_3. \label{tau}
\end{equation}
\begin{rem}
These relations \eqref{CrZW}, \eqref{CrWV}, \eqref{CrVZ} between $\{V, W, Z\}$ have the same form as those between $\{V, X, Z\}$ that define $m\mathfrak{H}$ and are given in \eqref{CZX}, \eqref{CXV}, \eqref{CVZ}. They actually have an extra parameter in front of $Z$ in the commutation relation between $V$ and $W$; this simply follows from the fact that the pencil parameter $\mu$ in $W$ reintroduces the free constant that we had taken away in the definition of $m\mathfrak{H}$ by using precisely this freedom (see \eqref{aut}) in the definition of $X$.

\end{rem}

\subsection{Rational Hahn algebra: $r\mathfrak{H} \hookrightarrow m\mathfrak{H}$}

Instead of using a linear combination of generators as in the preceding subsection, the embedding of the rational Hahn algebra in the meta one is achieved by supplementing $X$ and $Z$ with the generator $Y$ constructed as the following product:
\begin{equation}
    Y=XV.
\end{equation}
We are then led to:
\begin{prop}
A central extension of the rational Hahn algebra is embedded in the meta-Hahn algebra  by taking $X, Y=XV, Z$ as generators. The coefficients $\xi_i, i=0,\dots, 4$ entering in the relations \eqref{CMZX}, \eqref{CMXY} and \eqref{CMYZ} of $r\mathfrak{H}$ are given as follows in terms of the parameters $\eta_0, \eta_1, \eta_3$ of $m\mathfrak{H}$:
\begin{align}
    &\xi_0=\frac{1}{2}(\eta_3 - Q), \quad \xi_1=\eta_1+1, \quad \xi_2=\eta_0 +\eta_3 +1, \\
    &\xi_3=\eta_0 +\eta_1 + \eta_3 +1, \quad \xi_4 = 2\eta_3 +1.
    \end{align}
\end{prop}
The proof is again straightforward and simply relies on the repeated use of the commutation relations \eqref{CZX}, \eqref{CXV} and \eqref{CVZ} of $m\mathfrak{H}$. Note that here also the Casimir element $Q$ given in \eqref{Q} must be called upon to rexpress terms such as $VZ$ and $VZ^2$ in the generators. This explains how $Q$ arises in $\xi_0$ and why it is a central extension of $r\mathfrak{H}$ that is embedded in $m\mathfrak{H}$.

\section{Representations of the meta-Hahn algebra $m\mathfrak{H}$}

This section will develop certain aspects of the representation theory of the meta-Hahn algebra. We shall leave to another occasion its systematic study and will here be content with bringing to light the features that are most relevant to the bispectrality issues. We shall be concerned with finite dimensional representations on a real vector space $\mathfrak{V}_N$ of dimension $N+1$ with scalar product denoted by $( , )$.

We shall introduce a number of bases:
\begin{enumerate}
    \item The basis $\{d_n, n=0,\dots, N\}$ for the solution space of the GEVP
    \begin{equation}
       Xd_n=\lambda_n Zd_n;\label{GEVP1}
    \end{equation}
   \item The basis $\{d^*_n, n=0, \dots, N \}$ associated to the adjoint GEVP
   \begin{equation}
         X^T d^*_n=\mu_nZ^Td^*_n; \label{GEVP2}
   \end{equation}
    \item The eigenbasis $\{e_n, n=0,\dots, N\} $ of $V$, solutions of the EVP
    \begin{equation}
       Ve_n=\nu_n e_n;  \label{EVP1}
    \end{equation}
  \item The eigenbasis $\{f_n, n=0,\dots, N\}$ of the linear pencil $X+\mu Z$, with
  \begin{equation}
     (X+\mu Z)f_n= \rho _n f_n. \label{EVP2}
  \end{equation}
\end{enumerate}
We shall also consider the adjoint bases $\{e^*_n\}$ and $\{f^*_n\}$ such that $V^Te^*_n=\nu e^*_n$ and $(X^T +\mu Z^T)f^*_n=\rho_n f^*_n$. We shall determine the various spectra and the shapes of the matrices representing the generators in these bases.

The sets of eigenvalues $\{\lambda_n\}$ and $\{\mu_n\}$ coincide and are given by the roots of the characteristic polynomial of degree $N+1$ given by the determinant $|X-\lambda Z|$.

It is readily seen that the sets $\{d_n\}$ and $\{\Tilde{d}_n=Z^Td^*_n\}$ form two biorthogonal families of vectors:
\begin{equation}
    (d_m, \Tilde{d}_n)=0, \quad m \neq n \quad \text{with} \quad \Tilde{d}_n=Z^Td^*_n. \label{biorth}
\end{equation}
Indeed we have:
\begin{equation}
    (Xd_m, d^*_m)=\lambda_m(Zd_m, d^*_n)=(d_m, X^Td^*_n)=\lambda_n(d_m, Z^Td^*_n)=\lambda_n(Zd_m, d^*_n);
\end{equation}
and hence $(\lambda_m-\lambda_n)(Zd_m, d^*_n)=0$ which implies \eqref{biorth} assuming $\lambda_m \neq \lambda_n$.
\subsection {The GEVP basis $\{d_n\}$} \label{4.1}
We shall now obtain features of the matrices representing $X$, $Z$ and $V$ in the basis given by the vectors satisfying the GEVP $Xd_n=\lambda _n Zd_n$.
\begin{prop}
The eigenvalues $\lambda_n$ of the GEVP $Xd_n=\lambda_n Zd_n$ are linear and of the form
\begin{equation}
    \lambda_n=\alpha -n, \quad n=0,\dots N \qquad \text{with} \quad \alpha \in \mathbb{R.}\label{eigGEVP}
\end{equation}
$Z$ is a raising operator that acts as follows in the basis $\{d_n\}$:
\begin{equation}
    Zd_n=a_n d_{n+1}-d_n, \quad Zd_N=-d_N.\label{Zd}
\end{equation}
where $a_n$ is a $n$-dependent factor that depends on the normalization of the basis vectors.
\end{prop}
Proof. With the help of the commutation relation \eqref{CZX}, it is straightforward to obtain the intertwining relation
\begin{equation}
    (X-\Tilde{\lambda}_nZ)(Z+1)=(Z+1)(X-\lambda_nZ) \quad \text{with} \quad \Tilde{\lambda}_n=\lambda_n - 1.\label{IR}
\end{equation}
It follows that with $d_n$ solution of \eqref{GEVP1}, $(Z+1)d_n$ will be solution of the same GEVP with eigenvalue $\Tilde{\lambda}_n=\lambda_n - 1$.
The eigenvalues are thus incremented by unit steps and hence depend linearly on $n$. We shall take $Z+1$ to be a raising operator; ($Z^T+1)$ will therefore be lowering. This leads to the formula \eqref{eigGEVP} for the eigenvalues and allows to take the action of $Z$ on $d_n$ as given by \eqref{Zd}. This is a rare example of a GEVP that can be solved algebraically.
\begin{prop}\label{propVd}
The generator $V$ is represented in the basis $\{d_n\}$ by a lower Hessenberg matrix, that is, a matrix with zero entries above the first superdiagonal.
\end{prop}
Proof. Write generically $Vd_n=\sum_{m=0}^{N} V^{(\it{d})}_{m,n}d_m$. We shall use the constraints that relation \eqref{CVZ} imposes on the representations. Recall \eqref{Zd}. Given that $X$ has the same matrix form as $Z$  in view of the GEVP \eqref{GEVP1}, acting on $d_n$ with the right-hand side of \eqref{CVZ} will give an expression like $\kappa_nd_{n+1}+\zeta_nd_n$ where $\kappa_n$ and $\zeta_n$ are some coefficients.  The left-hand side will bring $[V,Z]d_n$ and the entries of the matrix $[V,Z]$ are given by
\begin{equation}
    [V,Z]^{(\it{d})}_{m,n} = a_nV^{(\it{d})}_{m,n+1} -a_{m-1}V^{(\it{d})}_{m-1,n}.
\end{equation}
For any matrix $A$, it is always assumed that $A_{m,n} =0$ whenever $m$ or $n$ are either smaller than $0$ or larger than $N$.

Consider first the entries above the diagonal. Since these are all zero on the right hand side, those of $[V, Z]$ must vanish. It is readily seen that this implies that all the entries of $V$ above the superdiagonal are $0$:
\begin{equation}
    V^{(\it{d})}_{n, n+k+1}=0, \qquad k=1, \dots, N-n-1. \label{upV}
\end{equation}
Indeed from $[V, Z]^{(\it{d})}_{1, k+1} =a_{k+1}V^{(\it{d})}_{1,k+2}=0$ for $k=1, \dots,N-2$, we find that \eqref{upV} is satisfied for $n=1$. Iterating over $n$ yields the desired result.

The entries on the diagonal and subdiagonal are subjected to more involved relations that would need to be solved were we to fully obtain the representation. Turn now attention to the entries below the subdiagonal. Again, since they have no match on the right of \eqref{CVZ}, the corresponding entries of $[V, Z]$ must vanish, which requires some of the matrix elements $V^{(\it{d})}_{m,n}$ of $V$ to satisfy
\begin{equation}
a_nV^{(\it{d})}_{m,n+1} -a_{m-1}V^{(\it{d})}_{m-1,n}=0, \qquad \text{for} \quad m > n+1, \quad m=n+2, \dots, N.\label{condbelow}
\end{equation}
This confirms that when acting on the vectors $\{d_n\}$, the only non-zero matrix elements of $V$ are above the superdiagonal.

The representations of the algebra $m\mathfrak{H}$ in the GEVP basis have another striking feature that we discuss next.
\begin{prop}\label{prop3}
In the basis $\{d_n\}$, the product $V^{(\it{d})}X^{(\it{d})}$ of the lower Hessenberg matrix $V^{(\it{d})}$ and the lower two-diagonal matrix $X^{(\it{d})}$ is tridiagonal.
\end{prop}
The proof of this proposition will make use this time of the relation \eqref{CXV} of $m\mathfrak{H}$. Recall that $Xd_n=(\alpha-n)Zd_n$, with $Zd_n$ given by \eqref{Zd}. It is straightforward to derive the following formulas for the matrix elements of $[X, V]$ and $\{V, Z\}$:
\begin{equation}
    [X, V]^{(\it{d})}_{mn}= (\alpha-m+1)a_{m-1}V^{(\it{d})}_{m-1,n}+(m-n)V^{(\it{d})}_{m,n}-(\alpha-n)V^{(\it{d})}_{m,n+1},
\end{equation}
\begin{equation}
    \{V, Z\}^{(\it{d})}_{m,n}=a_nV^{(\it{d})}_{m,n+1}-2V^{(\it{d})}_{m,n}+a_{m-1}V^{(\it{d})}_{m-1,n}.\label{VZanti}
\end{equation}
When considering the entries outside the diagonal and the subdiagonal when the left and the right hand sides of \eqref{CXV} are applied to $d_n$, only the terms $[X, V]$, $\{V, Z\}$ and $V$ contribute and lead to the following relations between  matrix elements of $V$:
\begin{equation}
    (\alpha-m)a_{m-1}V^{(\it{d})}_{m-1,n}+(m-n+1)V^{(\it{d})}_{m,n}-(\alpha-n+1)a_nV^{(\it{d})}_{m,n+1}=0, \quad m\neq\{n, n+1\}.\label{condition}
\end{equation}
Consider now the matrix elements of the product $VX$. It is easily seen that
\begin{equation}
    (VX)^{(\it{d})}_{m,n}=(\alpha-n)(a_nV^{(\it{d})}_{m,n+1}-V^{(\it{d})}_{m,n}).\label{VXmn}
\end{equation}
Owing to \eqref{upV}, above the superdiagonal, that is for $m<n-1$, we directly see that $(VX)_{m,n}=0$. For the entries below the subdiagonal, that is when $m>n+1$, we first observe that the conditions \eqref{condbelow} which we obtained from relation \eqref{CVZ} do apply. With their help, we see that equation \eqref{condition} which is a consequence of the defining relation \eqref{CXV} can be transformed into
\begin{equation}
    (m-n+1)(a_nV^{(\it{d})}_{m,n+1}-V^{(\it{d})}_{m,n})=0,
\end{equation}
as long as $m>n+1$.
Comparing with the expression \eqref{VXmn} for $(VX)^{(\it{d})}_{m,n}$, this allows to conclude that the matrix elements of $VX$ in the basis $\{d_n\}$ are also zero below the subdiagonal and therefore that remarkably VX is tridiagonal in the GEVP basis.

\subsection {Representation in the GEVP basis$\{d^*_n\}$}
This basis is defined by $X^Td^*_n=(\alpha-n)Z^Td^*_n$. The transpose of \eqref{IR} shows that
\begin{equation}
    Z^Td^*_n=-d^*_n+a^{(*)}_nd^*_{n-1} \label{ZT}
\end{equation}
where $a^{(*)}_n$ are some other $n$-dependent constants that are again dependent on the normalization of the basis vectors. This means that the only non-zero matrix elements of $Z$ in the basis $\{d^*_n\}$ are
\begin{equation}
    Z^{(*)}_{n+1,n}=a^{(*)}_n, \qquad Z^{(*)}_{n,n}=-1.
\end{equation}
Let $D$ be the diagonal operator
\begin{equation}
    Dd^*_n=(\alpha-n)d^*_n.
\end{equation}
(Note that we also have $Dd_n=(\alpha-n)d_n$.) In the basis $\{d^*_n\}$, we have $X^Td^*_n=Z^TDd^*_n$. This implies that the relation between $X$ and $Z$ in that basis is $X=DZ$ and we thus have:
\begin{prop} \label{propXd*}
The action of $X$ on the vectors $d^*_n$ is:
\begin{align}
    Xd^*_n&=DZd^*_n \label{X*}\\\nonumber \label{X*}
    &=X^{(\it{d}^*)}_{n+1,n}d^*_{n+1}+X^{(\it{d}^*)}_{n,n}d^*_n\\\nonumber
    &=(\alpha-n-1)a^{(*)}_nd^*_{n+1}-(\alpha-n)d^*_n.
    \end{align}
\end{prop}
The proof is immediate.

The representation of $V$ has the property:
\begin{prop}\label{4.5}
The action of $V$ in the basis $\{d^*_n\}$, $Vd^*_n=\sum_{m=0}^{N} V^{(\it{d}^*)}_{m,n}d^*_m$, is given by a lower Hessenberg matrix $V^{(\it{d}^*)}$.
\end{prop}
The proof follows the same lines as that of Proposition \ref{propVd}. The defining relation \eqref{CXV} imposes
\begin{equation}
     V^{(\it{d}^*)}_{n, n+k+1}=0, \qquad k=1, \dots, N-n-1, \label{upV*}
\end{equation}
and the conditions
\begin{equation}
   a^{(*)}_nV^{(\it{d}^*)}_{m,n+1} -a^{(*)}_{m-1}V^{(\it{d}^*)}_{m-1,n}=0, \qquad \text{for} \quad m > n+1, \quad m=n+2, \dots, N.\label{condbelow*}
\end{equation}
In addition, we find that in this basis, the analog of Proposition \ref{prop3} has the order of the operators $X$ and $V$ exchanged.
\begin{prop}
In the basis $\{d^*_n\}$, the product $X^{(\it{d}^*)} V^{(\it{d}^*)}$ is tridiagonal.
\end{prop}
The proof is similar to that of Proposition \ref{prop3}. We have
\begin{equation}
    [X, V]^{(\it{d}^*)}_{m,n}= (\alpha-m)a^{(*)}_{m-1}V^{(\it{d}^*)}_{m-1,n}+(m-n)V^{(\it{d}^*)}_{m,n}-(\alpha-n-1)V^{(\it{d}^*)}_{m,n+1}.
\end{equation}
The matrix elements $\{V,Z\}^{(\it{d}^*)}_{m,n}$ take the same form as \eqref{VZanti} with $V^{(\it{d})}$ replaced by $V^{(\it{d}^*)}$. For $m \neq \{n, n+1\}$, following \eqref{CXV},
\begin{align}
    &[X, V]^{(\it{d}^*)}_{m,n} - \{V,Z\}^{(\it{d}^*)}_{m,n} - V^{(\it{d}^*)}_{m,n} =\label{cond*}\\ \nonumber
    &  (\alpha-m-1)a^{(*)}_{m-1}V^{(\it{d^*})}_{m-1,n}+(m-n+1)V^{(\it{d}^*)}_{m,n}-(\alpha-n)a^{(*)}_nV^{(\it{d}^*)}_{m,n+1}=0.
\end{align}
Now observe that
\begin{equation}
     (XV)^{(\it{d}^*)}_{m,n}=(\alpha-n)(a_{m-1}V^{(\it{d}^*)}_{m-1,n}-V^{(\it{d}^*)}_{m,n}).
\end{equation}
Above the superdiagonal, for $m<n+1$, this expression vanishes in view of \eqref{upV*}. Below the subdiagonal, for $m>n+1$, with the help of \eqref{condbelow*}, condition \eqref{cond*} becomes
\begin{equation}
     (n-m-1)(a^{(*)}_{m-1}V^{(\it{d}^*)}_{m-1,n}-V^{(\it{d}^*)}_{m,n})=0,
\end{equation}
from where we see that the entries of $(XV)^{(\it{d}^*)}$ below the subdiagonal are also zero, thereby showing that $XV$ is tridiagonal in the basis $\{d^*_n\}$.

\subsection {The eigenbasis of $V$} \label{4.3}
We consider now the eigenbasis $\{e_n\}$ defined by \eqref{EVP1}. Let $Xe_n=\sum_{m=0}^{N}X^{(\it{e})}_{mn} e_m$ and $Ze_n=\sum_{m=0}^{N}Z^{(\it{e})}_{mn} e_m$.
Relevant observations regarding the representations of the meta-Hahn algebra in this basis are summarized here:
\begin{prop}
The eigenvalues $\nu_n$ of $V$ are quadratic and given by
\begin{equation}
    \nu_n=-n(n-\eta _1 -1).\label{eig}
\end{equation}
Moreover, $X$ and $Z$ act tridiagonally in the eigenbasis $\{e_n\}$ of $V$:
\begin{align}
    &Xe_n=X^{(\it{e})}_{n+1,n}e_{n+1}+X^{(\it{e})}_{n,n} e_n+X^{(\it{e})}_{n-1,n}e_{n-1}, \quad &X^{(\it{e})}_{-1,0}=0, \quad X^{(\it{e})}_{N+1,N}&=0;\label{3d1}\\
    &Ze_n=Z^{(\it{e})}_{n+1,n}e_{n+1}+Z^{(\it{e})}_{n,n} e_n+Z^{(\it{e})}_{n-1,n}e_{n-1}, \quad &Z^{(\it{e})}_{-1,0}=0, \quad Z^{(\it{e})}_{N+1,N}&=0.\label{3d2}
\end{align}
\end{prop}
Proof. The relations \eqref{CXV} and \eqref{CVZ} of $m\mathfrak{H}$ respectively impose the following conditions on the matrix elements of $X^{(\it{e})}$ and
$Z^{(\it{e})}$:
\begin{align}
    &(\nu_n -\nu_m -\eta_1)X^{(\it{e})}_{m,n}=(\rho_n +\nu_m -\eta_1)Z^{(\it{e})}_{m,n}+(\nu_n +\eta_0)\delta_{mn},\label{cond1}\\
    &(\nu_m +\nu_n -\eta_1)Z^{(\it{e})}_{m,n}=2X^{(\it{e})}_{m,n} + \eta_3 \delta_{mn}.\label{cond2}
\end{align}
Assume $m \neq n$. Expressing $X^{(\it{e})}_{m,n}$ in terms of $Z^{(\it{e})}_{m,n}$, using \eqref{cond2} and substituting in \eqref{cond1}, we see that the compatibility of these two equations requires that:
\begin{equation}
    (\nu_m -\nu_n )^2 + 2(\nu _m +\nu _n ) -\eta _1(\eta _1 +2)=0, \qquad m\neq n. \label{eigeq}
\end{equation}
For a given $n$, this quadratic relations allows for two possible values of $m$ which we will take to be $m=n \pm 1$ (without loss of generality). It follows that the action of $X$ and $Z$ will be tridiagonal as per \eqref{3d1} and \eqref{3d2}. Upon solving \eqref{eigeq} we readily see that the eigenvalues $\nu _n$ are given by \eqref{eig} (the second solution amounts to translating $n$ by $\eta_1 +1$).

The representations of $m\mathfrak{H}$ deserve a detailed study that will be carried out elsewhere. On their own, the representations of the subalgebra consisting of the deformed Jordan plane $\mathcal{J}_1$ generated by $X$ and $Z$ are already bound to be quite rich. It is striking that when the generator $V$ is added to enlarge $\mathcal{J}_1$ to form $m\mathfrak{H}$, the diagonalization of $V$ \textit{leads} to representations of $\mathcal{J}_1$ by tridiadiagonal matrices akin to those examined in \cite{tsujimoto2017tridiagonal} in the case of the $q$-oscillator algebra. One such representation of $\mathcal{J}_1$ has already been obtained in \cite{tsujimoto2020algebraic}.
We shall now view $\eta_0$ and $\eta_1$ as central charges and take them to have in this representation the following expressions
in terms of two parameters $\alpha, \beta$ and $N$ (where $N+1$ is the dimension):
\begin{equation}
   \eta_0 = (N-1-\beta)\alpha + \beta +1, \qquad \eta_1=N-1-\beta;\label{sv}
\end{equation}
we shall also parametrize $\eta_3$ in terms of $\alpha$ and $\beta$ by
\begin{equation}
   \eta_3=2\alpha - \beta -1. \label{sv1}
\end{equation}
It can be checked directly that the relations \eqref{CZX}, \eqref{CXV} and \eqref{CVZ} of $m\mathfrak{H}$ are satisfied when V is diagonal,
\begin{equation}
   V_{n,n}=-n(n-N+\beta),
\end{equation}
and $X$ and $Z$ are represented by $(N+1)\times(N+1)$ tridiagonal matrices with non-zero entries:
\begin{equation}
    X_{n+1,n}=\frac{n(\beta+n+1)(N-\beta-n)}{(N-\beta-2n)(N-\beta-2n-1)},  \label{X+}
\end{equation}
\begin{equation}
    X_{n,n}=-\alpha-n+\frac{\beta n(n-1)+{n}^2 (n-1)}{N-\beta-2n+1}-\frac {\beta n(n+1) +n(n+1)^2}{N-\beta-2n-1}, \label{Xd}
\end{equation}
\begin{equation}
   X_{n-1,n}= \frac {n(N-\beta-n)(N-n+1)}{(N-\beta-2n)(N-\beta-2n+1)},  \label{X-}
\end{equation}
 \begin{equation}
  Z_{n+1,n} = -\frac {(\beta+n+1)(N-\beta-n)}{(N-\beta-2n)(N-\beta-2n-1)},  \label{Z+}
 \end{equation}
 \begin{equation}
  Z_{n,n}=\frac {(n+1)\beta + (n+1)^2}{N-\beta-2n-1}-\frac{\beta n+n^2}{N-\beta-2n+1},  \label{ZD}
 \end{equation}
 \begin{equation}
   Z_{n-1,n}=-\frac {n( N-n+1) }{(N-\beta-2n) ( N-\beta-2n+1) }.   \label{Z-}
 \end{equation}
 The range of $n$ is the integers $\{0, 1, \dots, N\}$ and it is understood that $X_{N+1,N} = Z_{N+1,N}=0$.
The Casimir $Q$ of $m\mathfrak{H}$ given by \eqref{Q} then becomes
\begin{equation}
    Q=-2\alpha(\alpha-\beta) +2\alpha - \beta -1.
\end{equation}
Quite remarkably this representation of the meta-Hahn algebra that can be inferred from the representation found in \cite{tsujimoto2020algebraic} of the rational Hahn algebra allows to connect the deformed Jordan plane or the subalgebra generated by $X,Z$ to the Hahn polynomials.
Recall \cite{koekoek2010hypergeometric} that the recurrence relation
\begin{equation}
    -x Q_n(x)=A_nQ_{n+1}(x)-(A_n +C_n)Q_n(x)+C_nQ_{n-1}(x) \label{recurrence}
\end{equation}
of the Hahn polynomials $Q_n(x, \hat{\alpha}, \hat{\beta)}$ whose explicit expressions are given in \eqref{HahnOP}, has for coefficients
\begin{equation}
A_n =\frac{(n+\hat{\alpha}+\hat{\beta}+1)(n+\hat{\alpha}+1)(N-n)}{(2n+\hat\alpha+\hat{\beta}+1)(2n+\hat\alpha+\hat{\beta}+2)},\\ \label{An}
\end{equation}
\begin{equation}
  C_n =\frac{n(n+\hat{\alpha}+\hat{\beta}+N+1)(n+\hat{\beta})}{(2n+\hat\alpha+\hat{\beta})(2n+\hat\alpha+\hat{\beta}+1)}. \label{Cn}
\end{equation}
\begin{prop}
The Jacobi matrix $(X+\mu Z)_{m,n}$ obtained from \eqref{X+} - \eqref{Z-} is diagonalized by the orthogonal polynomials $\Tilde{Q}_n(y) = \sigma_n Q_n(x, \hat{\alpha}, \hat{\beta}, N)$ where $Q_n(x)$ are the Hahn polynomials and
\begin{align}
    x=&y-\alpha-\mu, \label{i} \\
    \hat{\alpha}=&-1-\mu, \label{ii} \\
    \hat{\beta}=&\beta + \mu -N, \label{iii}
\end{align}
and
\begin{equation}
    \sigma_n=\frac{(-1)^n (-N)_n}{(\hat{\alpha}+\hat{\beta}+N+1)_n}. \label{iv}
\end{equation}

\end{prop}
The proof is straighforward. Allowing for the normalization factors $\sigma_n$, the following relations are readily verified using formulas \eqref{i} to \eqref{iv}:
\begin{align}
    &(X_{n+1,n} + \mu Z_{n+1,n})\frac{\sigma_{n+1}}{\sigma_n}=-A_n, \\
    &(X_{n,n} + \mu Z_{n,n})=-\alpha -\mu +(A_n+C_n), \\
    &(X_{n-1,n} + \mu Z_{n-1,n})\frac{\sigma_{n-1}}{\sigma_n}=-C_n,
\end{align}
with $A_n$ and $C_n$ given by \eqref{An} and \eqref{Cn}. The conclusion of the proposition then follows from comparison with the recurrence relation \eqref{recurrence}.
This connection between the deformed Jordan plane and the Hahn polynomials takes us to the next subsection.

\subsection {The eigenbasis of the linear pencil $X+\mu Z$}
We shall now briefly discuss how the generators $X$, $Z$ and $V$ act on the vectors $f_n$ that satisfy $(X+\mu Z)f_n= \rho _n f_n$.
We have shown in Section 3 that the pencil $W=X+\mu Z$ together with $V$ defines an embedding of the Hahn algebra into $m\mathfrak{H}$. From the general theory \cite{terwilliger2004two}, it follows that when represented on the finite-dimensional space $\mathfrak{V}_N$, $V$ and $W$ form a Leonard pair. As a result, $V$ will act tridiagonally on the eigenvectors $f_n$ of $W$. We already know that $X$ and $Z$ separately act in a tridiagonal fashion in the eigenbasis of $V$; this will a fortiori be true for the pencil $W$ as should be from the Leonard perspective.

From the observation made in the last subsection that $W$ can be represented (essentially) by the Jacobi matrix of the Hahn polynomials we know that its spectrum will be linear. This can be found directly from the algebra \cite{granovskii1991quadratic}:
\begin{prop}
The spectrum of $W$ is given by
\begin{equation}
    \rho_n = an+b \label{rho}
\end{equation}
with $a$ and $b$ as affine parameters.
\end{prop}
The proof goes like this. We have $Wf_n = \rho _n f_n$ and we know that $V$ acts tridiagonally: $Vf_n = V^{(\it{f})}_{n+1,n}f_{n+1}+V^{(\it{f})}_{n,n}f_{n}+V^{(\it{f})}_{n-1,n}f_{n-1}$. $W$ and $V$ satisfy the relations of the Hahn algebra \eqref{121} and \eqref{212}. The last one reads:$[V,[W,V]]=a\{V,W\}+bV+c_2W+d_2\mathcal{I}$. Acting with both sides on $f_n$, the coefficient of $f_{n+2}$ yields the condition:
\begin{equation}
    \rho_{n+2} - 2\rho_{n+1} + \rho_n = 0,
\end{equation}
which has \eqref{rho} as solution.
\begin{rem}\label{rem4.1}
The parameters $a$ and $b$ can obviously be modified by an affine transformation of $W$. In the representation given in \eqref{X+} - \eqref{Z-}, we have as observed
$a=1$ and $b=-\alpha -\mu$.
\end{rem}
\begin{rem} \label{rem4.2}
It should be noted that the commutation relations of $m\mathfrak{H}$ show that this algebra does not admit representations where the three generators $X, Z, V$ are Hermitian (a fortiori symmetric) matrices. This is readily seen from the contradictions that follow by trying to impose such a constraint.
It is manifest however that the Hahn algebra $\mathfrak{H}$ defined by the relations \eqref{121}, \eqref{212} which are verified by $W$ and $V$ is compatible with Hermiticity. Under proper positivity conditions, it is known that $\mathfrak{H}$ admits unitary, in fact orthogonal, finite-dimensional representations. In these circumstances $W$ and $V$ can be symmetrized by diagonal similarity transformations and there is then no distinctions between the EVP and the adjoint EVP bases.
\end{rem}
\begin{rem}
The action of $Z$ (and of $X$) on the eigenvectors $\{f_n\}$ of $W$ can be obtained from those of $W$ and $V$ from the relations \eqref{CrZW}, \eqref{CrWV} and \eqref{CrWV} of the (non-standardized) meta-Hahn algebra that these operators verify together with $Z$. These actions are non-local, that is they typically correspond to full matrices.
\end{rem}
This is readily seen. From \eqref{CrWV} we have
\begin{equation}
  \{V,Z\} + \tau_2 Z = [W,V] - \tau_1 W - V - \tau_0 \mathcal{I}.
\end{equation}
Using \eqref{CrVZ}, we may transform this equation into
\begin{equation}
  (2V-\tau_2-\tau_1)Z = [W,V] +(2-\tau_1) W - V+(\tau_3 -\tau_0) \mathcal{I}. \label{Heun}
\end{equation}
Note that the right-hand side (rhs) of \eqref{Heun} is an algebraic Heun operator \cite{grunbaum2018algebraic} and is tridiagonal in both the eigenbases of $V$ and $W$ \cite{nomura2007linear}. Clearly, $Z$ can then be obtained by multiplying this operator on the rhs by the inverse of the factor of $Z$ on the left, that is by $(2V-\tau_2-\tau_1)^{-1}$. The action of $X$ is then given via $X=W-\mu Z$.

\section{Bispectrality}

We are now ready to discuss the bispectral properties of the functions defined by overlaps of the form $(g_m,h_n), m,n=0,\dots,N$ where $g_m$ and $h_m$ are some of the basis vectors in $\mathfrak{V}_N$ considered in the last section. Recall that the scalar product is taken to be real; to be clear, it will be defined as follows using row and column vectors:
\begin{equation}
    (g_m,h_n)=g^T_mh_n=h^T_ng_m=(h_n,g_m).
\end{equation}
The transpose of an operator $R$ on $\mathcal{V}_N$ is such that
\begin{equation}
    (g_m, Rh_n) = (R^Tg_m,h_n)
\end{equation}
which is consistent with $g^T_mRh_n =(R^Tg_m)^Th_n$. The overlap $(g_m,h_n)$ involves two discrete indices $m$ and $n$ that will play the role of degree and variable (or vice-versa for the dual functions).

A key result obtained in \cite{zhedanov1999biorthogonal} states that solutions to GEVPs involving two tridiagonal matrices are biorthogonal rational functions. This is somehow the analog to Favard's theorem \cite{chihara2011introduction} connecting EVPs for Jacobi matrices to orthogonal polynomials. The overlaps stemming from the representations of the meta-Hahn algebra will lead to functions of those two kinds. The issue will be to recognize first that these characterizations apply and second, that the algebra $m\mathfrak{H}$ entails bispectrality. While more space will be dedicated to the less understood biorthogonal rational functions, we wish to stress that the meta-Hahn algebra underpins in a unified way the bispectral properties of these two classes of functions.

\subsection {Biorthogonal rational functions}

Consider the functions defined as follows as overlaps between the GEVP bases and the eigenbasis of $V$ (and $V^T$):
\begin{equation}
    U_m(n)=(e_m, d^*_n) \qquad \text{and} \qquad \Tilde{U}_m(n)=(e^*_m, Z d_n).
\end{equation}
where the vectors $e^*_n$ are the eigenvectors of the transpose of $V$: $V^Te^*_n=\nu_n e^*_n$.
Before we show that these are rational functions and identify their bispectral properties, let us first observe that they are biorthogonal.
\subsubsection{Biorthogonality}\label{sub}
We have
\begin{equation}
(d_n, Z^Td^*_k)=(Zd_n,d^*_k)=w^{-1}_n \delta _{nk},  \qquad k,n=0, \dots,N \label{normd}
\end{equation}
and hence
\begin{equation}
\sum_{n=0}^{N} (u,Zd_n)(d^*_n, v)w_n=(u,v).
\end{equation}
It thus follows that
\begin{equation}
\sum_{n=0}^{N}\Tilde{U}_k(n)U_m(n)w_m=\sum_{n=0}^{N}(e^*_k, Zd_n)(d^*_n, e_m)w_n=(e^*_k, e_m)=0, \quad k\neq m. \label{utildeu}
\end{equation}

\subsubsection{Bispectral equations for $U_m(n)$}

From the definition of $U_m(n)$ and the fact that the vectors $d^*_n$ satisfy $X^Td^*_n=\lambda_n Z^Td^*_n$, we have
\begin{equation}
    (e_m, (X^T -\lambda_n Z^T)d^*_n)=0
\end{equation}
from where it follows that
\begin{equation}
  ((X-\lambda_n Z)e_m, d^*_n)=0.
\end{equation}
Recalling that $X$ and $Z$ have tridiagonal actions on the eigenvectors of $V$ and that the GEVP eigenvalue is $\lambda_n=\alpha-n$, we find the following result:
\begin{prop}
The function $U_m(n)$ satisfies the following recurrence relation of GEVP form
\begin{align}
  &X^{(\it{e})}_{m+1,m}U_{m+1}(n)+X^{(\it{e})}_{m,m}U_m(n)+X^{(\it{e})}_{m-1,m}U_{m-1}(n)=\nonumber \\
  &(a-n)\left[Z^{(\it{e})}_{m+1,m}U_{m+1}(n)+Z^{(\it{e})}_{m,m}U_m(n)+Z^{(\it{e})}_{m-1,m}U_{m-1}(n)\right].\label{recur}
\end{align}
\end{prop}

Let us now identify the other equation of the bispectral system. In a similar fashion we have
\begin{equation}
    ((V-\nu_m)e_m, d^*_n)=0
\end{equation}
since $e_n$ is an eigenvector of $V$. Multiplying the first factor by $X$ allows to write the identity $(X(V-\nu_m)e_m, d^*_n)=0$ which leads to
\begin{equation}
    (e_m, (V^TX^T-\nu_mX^T)d^*_n)=0.
\end{equation}
Since $Y=XV$ is tridiagonal in the GEVP basis $\{d^*_n\}$, the same will be true for the transposed matrix $Y^T=V^TX^T$. We recall that $\nu_m=-m(m-\eta_1 -1)$ and that (see \ref{propXd*})
\begin{equation}
    X^Td^*_n=(\alpha-n)[-d^*_n+a^{(*)}_n d^*_{n-1}].
\end{equation}
This leads to to the following observations:
\begin{prop}
The function $U_m(n)=(e_m, d^*_n)$ is a biorthogonal rational function that satisfies the difference equation of GEVP type:
\begin{align}
    &Y_{n,n+1}U_m(n+1)+Y_{n,n}U_m(n)+Y_{n,n-1}U_m(n-1)=\nonumber\\
    &-m(m-\eta_1 -1)(\alpha -n)\left[(\alpha -n)(-U_m(n)+a^{(*)}_nU_m(n-1))\right].\label{diff}
\end{align}
\end{prop}
It was already proved to be biorthogonal. That it is a rational function results from the fact that it verifies the above difference GEVP involving two tridiagonal matrices (one being in fact two-diagonal).
These observations are summarized in:
\begin{prop}
The rational function $U_m(n)$ is bispectral since it satisfies the recurrence relation \eqref{recur} and the difference equation \eqref{diff}. These properties are algebraically encoded in the meta-Hahn algebra.
\end{prop}

\subsubsection{The bispectral equations of $\Tilde{U}_m(n)$}
We may describe in the same spirit the bispectrality of the functions $\Tilde{U}_m(n)=(e^*_m, Zd_n)$. First observe that the vectors
 $\Tilde{d}_n=Zd_n$ satisfy the GEVP
\begin{equation}
    (\Tilde{X}-\lambda_n\Tilde{Z}){d}_n=0,
\end{equation}
with
\begin{equation}
    \Tilde{X}=X+Z+1, \qquad \Tilde{Z}=Z.
\end{equation}
This simply follows from $[Z,X]=Z^2+Z$. Let $\Tilde{V}=V$. The generators $\Tilde{X}$, $\Tilde{Z}$  and  $\Tilde{V}$ verify the relations \eqref{CrZW}, \eqref{CrWV}, \eqref{CrVZ} with $W$, $Z$ and $V$ replaced respectively by $\Tilde{X}$, $\Tilde{Z}$  and  $\Tilde{V}$ and the parameters $\tau$ given by (differently from \eqref{tau}): $\tau_1=2\eta_1 -2$, $\tau_2=2-3\eta_1$, $\tau_0=2+\eta_0 -\eta_1 -\eta_3$, $\tau_3=\eta_3 -2$. So apart from the fact that the parameters are different, the structure of the relations obeyed by $\Tilde{X}$, $\Tilde{Z}$  and  $\Tilde{V}$ is the same as that of the relations between $X$, $Z$ and $V$. We observe that in the new notation:
\begin{equation}
    \Tilde{U}_m(n)=(\Tilde{e}^*_n,\Tilde{d}_n) \qquad \text{with} \qquad \Tilde{V}^T\Tilde{e}^*_n=\nu_n\Tilde{e}^*_n.
\end{equation}
The essential features of the representations in the bases $\{\Tilde{d}_n\}$ and $\{\Tilde{e}^*_n\}$ of the algebra generated by  $\Tilde{X}$, $\Tilde{Z}$  and  $\Tilde{V}$ will hence be the same as those of the meta-Hahn algebra in the bases $\{d_n\}$ and $\{e^*_n\}$; namely, $\Tilde{X}$ and $\Tilde{Z}$ will be lower bidiagonal in the basis $\{\Tilde{d}_n\}$ and tridiagonal in the basis $\{e^*_n\}$, $V$ will be lower Hessenberg and $VX$ will be tridiagonal in the basis $\{\Tilde{d}_n\}$. This is confirmed by going again through the arguments that led to the propositions in Subsection \ref{4.1} and \ref{4.3}. Indeed, these only rely on the structure of the commutation relations which has not been modified by going to the tilded generators. Of course, the actual entries of the matrices representing the different operators $\Tilde{X}$, $\Tilde{Z}$  and  $\Tilde{V}$ will be affected by the changes in the parameters.

These observations make clear that the biorthogonal partners $\Tilde{U}_m(n)$ are also as expected, bispectral rational functions whose properties are determined by an algebra which is essentially $m\mathfrak{H}$. As a result,
\begin{equation}
    ((\Tilde{X}^T-\lambda_n\Tilde{Z}^T)\Tilde{e}^*_m, \Tilde{d}_n)=0
\end{equation}
leads to their recurrence relation while
\begin{equation}
    (\Tilde{e}^*_m,(\Tilde{V}\Tilde{X}-\Tilde{\nu}_m{X})\Tilde{d}_n)=0
\end{equation}
amounts to their difference equation.

\subsection{Bispectral orthogonal polynomials} \label{5.2}
The last functions we wish to consider are given by the overlaps
\begin{equation}
 S_m(n)=(f^*_n, e_m) \quad\text{and} \quad \Tilde{S}_m(n)=(f_n, e^*_m)
\end{equation}
between the eigenvectors of $W^T=X^T+\mu Z^T$ and those of $V$ and between reciprocally the eigenvectors of $W=X+\mu Z$ and of $V^T$. With
\begin{equation}
    (e^*_m, e_n)=\kappa ^{-1}_n \delta _{m,n}, \qquad  (f^*_m,f _n)=\zeta ^{-1}_n \delta _{m,n},
\end{equation}
and
\begin{equation}
    \sum_{n=0}^N (u, f_n)(f^*_n,v)\zeta_n=(u,v), \qquad  \sum_{m=0}^N (u, e_m)(e^*_m,v)\eta_n=(u,v),
\end{equation}
we have the relations
\begin{align}
    &\sum_{n=0}^N \Tilde{S}_m(n)S_k(n)\zeta_n=\sum_{n=0}^N (f_n, e^*_m)(f^*_n, e_k)\zeta _n = \kappa^{-1}_m \delta_{m,k}, \label{orth1}\\
    &\sum_{m=0}^N \Tilde{S}_m(n)S_m(k)\kappa_n=\sum_{m=0}^N (f_n, e^*_m)(f^*_k, e_m)\kappa_m = \zeta^{-1}_n \delta_{n,k}.\label{orth2}
    \end{align}
    Given that $V$ and $W$ verify the commutation relations of the Hahn algebra, within the appropriate range of the parameters, we expect the functions $S_m(n)=(f^*_n, e_m)$ and $\Tilde{S}_m(n)=(f_n, e^*_m)$ to be both related to the same (dual) Hahn polynomials so that \eqref{orth1} and \eqref{orth2} are orthogonality relations for a single family of polynomials. To convince oneself that $S_m(n)=(f^*_n, e_m)$ and $\Tilde{S}_m(n)=(f_n, e^*_m)$ can be proportional to these polynomials, one shall recall Remark \ref{rem4.2}
    as a reminder that the Hahn algebra generators can be symmetrized thus removing the difference between the starred and the non-starred bases and rendering $S_m(n)$ and $\Tilde{S}_m(n)$ identical. We would nevertheless like for generality and practical reasons to keep the adjoint bases and still conclude that $S_m(n)$ and $\Tilde{S}_m(n)$ are expressed in terms of the same polynomials. An argument to that effect goes like this.
From $\rho_n(f^*_n, e_m)=(W^Tf^*_n, e_m)=(f^*_n, We_m)$ we have
   \begin{equation}
       \rho_n S_m(n)=W^{(\it{e})}_{m+1,m}S_{m+1}(n)+W^{(\it{e})}_{m,m}S_{m}(n)+W^{(\it{e})}_{m-1,m}S_{m-1}(n).\label{recS}
   \end{equation}
   We may thus write
   \begin{equation}
       S_m(n)=S_0(n)P_m(n) \quad \text{with} \quad P_0(n)=1, \label{fact}
   \end{equation}
   and conclude from Favard's theorem that since they obey a three-term recurrence relation, the functions $P_m(n)$ are orthogonal polynomials if
   \begin{equation}
      W^{(\it{e})}_{m+1,m} W^{(\it{e})}_{m-1,m} > 0
   \end{equation}
   and the matrix elements of $W$ are real. According to \eqref{orth1}, the functions $\Tilde{S}_m(n)$ are themselves orthogonal (with weight factors) to these polynomials and we can hence conclude that we must have similarly to  \eqref{fact} $\Tilde{S}_m(n)=h(m)\Tilde{S}_0(n)P_m(n)$ since two systems $\{P_m\}$, $\{Q_m\}$, $m=0,\dots,N$ of orthogonal polynomials that are mutually orthogonal must be related by $P_m=h(m)Q_m$ with $h(m)$ a non-zero function of $m$. This is the desired conclusion, namely, the core of the functions $S_m(n)$ and $\Tilde{S}_m(n)$ is a single family of orthogonal polynomials. The identity $((W-\rho_n)f_n, e^*_m)=(f_n, (W^T-\rho_n)e^*_m)=0$ leads to
   \begin{equation}
      \rho_n \Tilde{S}_m(n)=W^{(\it{e}^*)T} _{m+1,m}\Tilde{S}_{m+1}(n)+W^{(\it{e}^*)T}  _{m,m}\Tilde{S}_m (n)+W^{(\it{e}^*)T} _{m-1,m}\Tilde{S}_{m-1}(n)
   \end{equation}
   which must be equivalent to \eqref{recS} and hence implies that $h(k)W^{(\it{e}^*)T} _{k,m}=W^{(\it{e})} _{k,m}$.

   The dual identities $(f_n^*,(V-\nu_m)e_m)=((V^T-\nu_m)f^*_n, e_m)=0$ and $(f_n,
(V^T-\nu_m)e^*_m)=((V-\nu_m)f_n, e^*_m)=0$ entail the bispectrality of $S_m(n)$ and
$\Tilde{S}_m(n)$ and the same difference equation for the polynomials $P_m(n)$ with the
implication that $S_0(n)V^{(\it{f}^*)T} _{k,n}=\Tilde{S}_0(m)V^{(\it{f})} _{k,n}$. The
equation for $S_m(n)$ is:
   \begin{equation}
     \nu_m S_{m}(n)=V^{(\it{f^*})T}_{n+1,n}  S_{m}(n+1)+V^{(\it{f^*})T}_{n,n} S_{m}(n)+V^{(\it{f^*})T}_{n-1,n}S_{m}(n-1). \label{diffS}
   \end{equation}
   We note that \eqref{orth1} and \eqref{orth2} embody the orthogonality relations of the polynomials $P_m(n)$ and their duals. These observations can be summed up as follows.
  \begin{prop}
   The functions $S_m(n)=(f^*_n, e_m)$ are bispectral and obey the three-term recurrence and difference relations \eqref{recS} and \eqref{diffS}; they are proportional to orthogonal polynomials.  The functions $\Tilde{S}_m(n)=(f_n, e^*_m)$ are also expressed in terms of the same polynomials.
   \end{prop}
   Stressing the following point might be appropriate as we conclude this section. To discuss special functions associated to GEVPs, the adjoint basis is required. This basis is different from that of the original problem. For consistency, we have also introduced the adjoint bases of the EVPs. A reader familiar with Leonard pairs knows that given such a set, say $(V,W)$, orthogonal polynomials arise through expanding the vectors of one eigenbasis in terms of the vectors of the other eigenbasis. Here the direct and the adjoint eigenbases can be related through symmetrizing factors. Expressions such as $(f^*_n, e_m)$ provide (up to the normalization factors) the expansion coefficients of the vector $f_n$ over the basis formed by the vectors $e_m, m=0,\dots, N$. Concrete examples of such computations will be provided in the next section.

\section{Models and the special functions associated to the meta-Hahn algebra} \label{6}

We shall now provide realizations of the elements of representation theory developed so far and in this way identify explicitly the special functions whose properties are rooted in $m\mathfrak{H}$.
\subsection {A differential realization}\label{6.1}
Consider the following operators (we shall use the same symbols for the abstract generators and their representations):
\begin{equation}
 Z=(x-1) \mathcal{I}, \label{Z_D}
\end{equation}
\begin{equation}
X = x(1-x)\frac{d}{dx} -\alpha  \mathcal{I},
\label{X_D}
\end{equation}
\begin{equation}
 V= x \left( 1-x \right)  \frac{d^2}{dx^2 }+\left[x \left( N-1-\beta \right) -N \right] \frac{d}{dx}  \label{V_D}.
\end{equation}
It is readily seen that they obey the commutation relations \eqref{CZX}, \eqref{CXV}, \eqref{CVZ} of $m\mathfrak{H}$ with the parameters assigned as in \eqref{sv} and \eqref{sv1}. Given this representation, the corresponding model for the basis $\{d_n\}$
is straightforwardly obtained. It is easily verified that the functions
\begin{equation}
    d_n(x;\alpha)=\gamma_n x^n(1-x)^{-\alpha}, \quad n=0, 1, \dots N, \label{dreal}
\end{equation}
with $\gamma_n$ a normalization factor, satisfy the differential GEVP  $Xd_n(x;\alpha)=(\alpha-n)Zd_n(x;\alpha)$.
 It is a matter of computation to obtain the actions of $Z$, $X$ and $V$ on these basis vectors. One readily sees that
\begin{equation}
    Zd_n(x;\alpha)=\frac{\gamma_n(\alpha)}{\gamma_{n+1}(\alpha)}d_{n+1}(x;\alpha)-d_n(x;\alpha)=-\frac{\gamma_n(\alpha)}{\gamma_n(\alpha +1)}d_n(x;\alpha +1)
\end{equation}
which is in line with \eqref{Zd} with $a_n=\frac{\gamma_n(\alpha)}{\gamma_{n+1}(\alpha)}$. One also finds
\begin{align}
   Xd_n(x;\alpha)&=(\alpha -n)\frac{\gamma_n(\alpha)}{\gamma_{n+1}(\alpha)}d_{n+1}(x;\alpha)+(n-\alpha)d_n(x;\alpha) \\
   &=(n-\alpha)\frac{\gamma_n(\alpha)}{\gamma_n(\alpha +1)}d_n(x;\alpha +1)
\end{align}
from where one confirms that $Xd_n(x;\alpha)=(\alpha-n)Zd_n(x;\alpha)$. We shall not record the action of $V$ on $d_n(x;\alpha)$ as we shall give below that of $V^T$ on $d^*_n(x;\alpha)$ which has similar features.

The adjoint bases will be connected with the formal Lagrange adjoints of the above operators and we will keep the practice of apposing an asterix to the corresponding solutions of the adjoint GEVPs or EVPs. We have
\begin{equation}
    Z^T=(x-1)\mathcal{I}, \label{Z_DT}
\end{equation}
\begin{equation}
    X^T=-x(1-x)\frac{d}{dx}+(2x-\alpha -1)\mathcal{I},\label{X_DT}
\end{equation}
\begin{equation}
    V^T=x(1-x)\frac{d^2}{dx^2}+\left[N+2-x(N+3-\beta)\right]\frac{d}{dx}-(N+\beta -1)\mathcal{I}.
\end{equation}
The solutions of $(X^T -(\alpha-n)Z^T)d^*_n(x;\alpha)=0$ are recognized to be:
\begin{equation}
    d^*_n(x;\alpha)=\gamma^*_n(\alpha)x^{-1-n}(1-x)^{-1+\alpha}, \quad n=0, 1\dots N,
\end{equation}
where $\gamma^*(\alpha)$ are the normalizing coefficients of these basis vectors.
The actions of $Z^T$ and $X^T$ on $ d^*_n(x;\alpha)$ are in this case also simple to derive and one gets:
\begin{equation}
    Z^T d^*_n(x;\alpha)=- d^*_n(x;\alpha)+\frac{\gamma^*_n(\alpha)}{\gamma^*_{n-1}(\alpha)}d^*_{n-1}(x;\alpha)=-\frac{\gamma^*_n(\alpha)}{\gamma^*_n(\alpha +1)}d^*_n(x;\alpha +1), \label{ZTd*}
\end{equation}
\begin{equation}
    X^Td^*_n(x;\alpha)=(n-\alpha)d^*_n(x;\alpha)+(\alpha-n)\frac{\gamma^*_n(\alpha)}{\gamma^*_{n-1}(\alpha)}d^*_{n-1}(x;\alpha)=(n-\alpha)\frac{\gamma^*_n(\alpha)}{\gamma^*_n(\alpha +1)}d^*_n(x;\alpha+1), \label{XTd*}
\end{equation}
which fits with \eqref{ZT} and \eqref{X*} and sets
\begin{equation}
    a^*_n=\frac{\gamma^*_n(\alpha)}{\gamma^*_{n-1}(\alpha)}. \label{normc}
\end{equation}
Finding the action of $V^T$ in this adjoint basis requires a little more work but is also straightforward and one gets:
\begin{align}
    &V^Td^*_n(x;\alpha)=-(N-n)(n+1)\frac{\gamma_n(\alpha)}{\gamma_{n+1}(\alpha)}d^*_{n+1}(x;\alpha)\nonumber\\
    &+\left[(N-n)(n-2\alpha+\beta+2)+N(\alpha-\beta-1)\right]d^*_n(x;\alpha)\nonumber\\
    &+(\alpha-\beta-1)(\alpha-1)\left[ \frac{\gamma_n}{\gamma_{n-1}}d^*_{n-1}(x;\alpha)
    +\frac{\gamma_n}{\gamma_{n-2}}d^*_{n-2}(x;\alpha)+ \dots + \frac{\gamma_n}{\gamma_0}d^*_0(x;\alpha) \right].\label{VTd*}
\end{align}
In conformity with Proposition \ref{4.5}, we see that $V^T$ is a upper Hessenberg matrix. We can trace the origin of the terms extending upwards beyond the diagonal to the presence of a factor $(1-x)^{-1}=1+x+x^2+\dots$ which yields the lowering sequence $d^*_n, d^*_{n-1}, \dots, d^*_0$ when it multiplies $d^*_n$. Note that truncation must be enforced when $d^*_0$ is reached.

The scalar product $(f(x), g(x))$ for $f(x), g(x) \in \mathfrak{V}_N$ is realized in this model by integrating over a closed contour in the complex plane. Let $\Gamma$ be a circle $|x|= a<1$. With the functions taken to be on $\mathbb{C}$, we define
\begin{equation}
     (f(x), g(x))=\frac{1}{2\pi i} \oint_{\Gamma} f(x)g(x)dx, \quad f, g \in \mathcal{V_N}.\label{sp}
\end{equation}
It is indeed readily seen using Cauchy's theorem that
\begin{equation}
    (d^*_m, Zd_n)= \gamma^*_m\gamma_n\frac{1}{2\pi i} \oint_{\Gamma} x^{(n-m-1)}dx=\gamma^*_m\gamma_m\; \delta_{m,n}. \label{d*d}
\end{equation}
We shall remark that we had introduced generically the normalization factor $w^{-1}_n$ in $(d^*_k, Zd_n)=w^{-1}_n \delta_{kn}$ (see \eqref{normd}), we thus have the identification
\begin{equation}
    w^{-1}_n=\gamma^*_n\gamma_n. \label{wgg}
\end{equation}

The bases $\{e_m\}$ and $\{e^*_m\}$ are realized as follows in this model. The solutions of the EVP $Ve_n(x; \beta, N)= \nu_n e_n(x; \beta, N)$, where $\nu_n=-n(n-N+\beta)$, are given by the polynomials:
\begin{equation}
   e_m(x;\beta, N) = \delta_m \; {_2}F_1 \left( {-m, m-N+\beta \atop  -N } ; x\right), \quad m=0,1,\dots, N. \label{ereal}
\end{equation}
Note that these are Jacobi polynomials where the parameter $\alpha$ is a negative integer and that the usual orthogonality relation is not valid since the corresponding integral is not defined. (See \cite{kuijlaars2005orthogonality} for a discussion of Jacobi polynomials with general parameters.) These functions are rather orthogonal under the scalar product given in \eqref{sp} to a set of Laurent polynomials that will form the adjoint basis $\{e^*_m\}$. Indeed
the EVP  $V^Te^*_m(x;\beta)=\nu_n e^*_m(x;\beta)$ admits the following solutions
\begin{equation}
    e^*_m(x;\beta, N)= \delta^*_m \; x^{-1-N}{_2}F_1 \left( {m-N, -m-\beta \atop  -N } ; x\right), \quad m=0,1,\dots, N, \label{e*real}
\end{equation}
where $\delta^*_m$ are normalization constants.
That they form the basis which is such that
\begin{align}
   (e^*_m, e_n)=& \delta^*_m\delta_n\frac{1}{2\pi i} \oint_{\Gamma}dx\;x^{-1-N}{_2}F_1 \left( {m-N, -m-\beta \atop  -N } ; x\right){_2}F_1 \left( {-n, n-N+\beta \atop  -N } ; x\right)\nonumber \\
 =&\xi_n \delta_{m,n}
\end{align}
with $\xi_n$ some constants is shown in Appendix \ref{App}.

\subsection{The rational functions of Hahn type}
We shall now identify explicitly the rational functions whose bispectrality is encoded in $m\mathfrak{H}$ and indicate how their characterization largely
follows from this algebraic framework.
\subsubsection{$U_m(n)$}
We shall use the differential realization to determine the functions $U_m(n)=(e_m, d^*_n)$. In the model described in Section \ref{6.1},
\begin{equation}
U_m(n)=\frac{1}{2\pi i} \oint_{\Gamma}dx\;\gamma^*_n\delta_m \;x^{-1-n}(1-x)^{-1+\alpha}{_2}F_1 \left( {-n, n-N+\beta \atop  -N } ; x\right).
\end{equation}
We hence need to find the residue or the coefficient of $x^{-1}$ in the Laurent series of the integrand or equivalently, the coefficient of $x^n$ in the expansion:
\begin{align}
(1-x)^{\alpha -1}  \: {_2}F_1 \left( {-n, n-N+\beta \atop  -N } ; x\right) =\frac{1}{\delta_m}\sum_{\ell} \frac{1}{\gamma^*_n}U_m(\ell) x^{\ell} \label{expand1}
\end{align}
which only involves non-negative powers and where we have formally extended the domain of $U_m$.
Using the binomial formula
\begin{equation}
 (1-x)^{\alpha} = \sum_{k=0}^{\infty} \frac{(-\alpha)_k x^k}{k!}, \label{binom}
\end{equation}
one readily obtains
\begin{equation}
 U_m(n) = \delta_m\gamma^*_n\frac{(1-\alpha)_n}{n!} \: {_3}F_2 \left(  { -n,-m, m-N+\beta \atop -N, \alpha-n }  ; 1 \right). \label{Umn}
\end{equation}
With $n$ as the variable (instead of $x$), this is seen to be, up to a normalization factor, the rational Hahn function $\mathcal{U}_m(n;\alpha,\beta,N)$ defined in \eqref{ratHahn} multiplied by $\delta_m\gamma^*_n\frac{(1-\alpha)_n}{n!}$
\subsubsection{A finite difference realization and contiguity relations}
A natural choice for the normalization constants for the vectors $d^*_n(x;\alpha)$ is
\begin{equation}
    \gamma^*_n=\frac{n!}{(1-\alpha)_n} \label{gamma}
\end{equation}
in which case $U_m(n) = \delta_m \: {_3}F_2 \left(  { -n,-m,\; m-N+\beta \atop -N,\; \alpha-n }  ; 1 \right)$. We then find
\begin{equation}
    \frac{\gamma^*(\alpha)}{\gamma^*_{n-k}(\alpha)}=\frac{n!}{(n-k)!}\frac{(1-\alpha)_{n-k}}{(1-\alpha)_n}=\frac{(-n)_k}{(\alpha-n)_k}
\end{equation}
with the help of the standard relations
\begin{equation}
    (-n)_k=(-1)^k\frac{n!}{(n-k)!} \quad \text{and} \quad (a)_{n-k}=(-1)^k \frac{(a)_n}{(1-n-a)_k}.\label{pochrel}
\end{equation}

Define the action of $Z$, $X$, and $V$ on $U_m(n)$ by
\begin{equation}
    \mathcal{O}U_m(n)=(e_m, \mathcal{O}^Td^*_n)
\end{equation}
for $\mathcal{O}$ any of $Z$, $X$, or $V$. Using the notation $T^{\pm}f(n)=f(n\pm1)$ and $\mathcal{I}f(n)=f(n)$, from \eqref{ZTd*}, \eqref{XTd*} and \eqref{VTd*}, we obtain the following realization of the meta-Hahn algebra generators in terms of difference operators:
\begin{equation}
    Z=-\mathcal{I}+\frac{n}{n-\alpha}T^-,
\end{equation}
\begin{equation}
    X=(n-\alpha)\mathcal{I}-nT^-,
\end{equation}
\begin{align}
    V=&-(N-n)(n-\alpha+1)T^+ +\left[(N-n)(n-2\alpha+\beta+2)+N(\alpha-\beta-1)\right]\mathcal{I} \nonumber \\
    &+(\alpha-\beta-1)(\alpha-1)\sum_{k=1}^N\frac{(-n)_k}{(\alpha -n)_k}(T^-)^k.
\end{align}
It is understood that
\begin{equation}
    (T^-)^k U_m(n)=0 \quad \text{for} \quad k>n.
\end{equation}
Note that in this model
\begin{equation}
    X=(\alpha-n)Z. \label{XZdiff}
\end{equation}
It is then further observed still in the framework of this realization, that
\begin{equation}
    Y=XV=A_1(n)T^++A_2(n)T^-+A_0(n), \label{Ydiff}
\end{equation}
where
\begin{align}
    A_1(n)=&(n-\alpha)(n-N)(n+1-\alpha), \nonumber \\
    A_2(n)=&n(n-\alpha)(n-\alpha +\beta -N), \nonumber \\
    A_0(n)=&(n-\alpha)(-2n^2+(2\alpha-1+2N-\beta)n-N(\alpha-1)).
\end{align}
This is the translation of the fact observed abstractly that the product $V^TX^T$ acts in a tridiagonal fashion in the basis $\{d^*_n\}$.
We thus recover the difference realization introduced in \cite{tsujimoto2020algebraic} as the starting point to the analysis of the bispectral properties of the rational Hahn functions \eqref{ratHahn}. The basis vectors which obey the GEVP $Xf(n)=(\alpha-n)Zf(n)$ in this difference realization are the rational Hahn functions themselves.

We can also recover contiguity relations for the rational functions $U_m(n;\alpha, \beta, N)$ using this framework. We have
\begin{equation}
    \frac{\gamma^*_n(\alpha)}{\gamma^*(\alpha +1)}=\frac{\alpha}{\alpha-n}
\end{equation}
and hence from the last equality in \eqref{XTd*} we see that
\begin{equation}
    X^Td^*_n(x,\alpha)=-\alpha d^*_n(x;\alpha +1).
\end{equation}
The relation
\begin{equation}
    XU_m(n;\alpha)=(e_m(x;\beta), X^Td^*_n(x;\alpha))=\frac{\alpha}{\alpha-n}(e_m(x,\alpha), d^*_n(x, \alpha +1))=\frac{\alpha}{\alpha -n}U_m(n;\alpha +1)
\end{equation}
gives the contiguity relation
\begin{equation}
    (n-\alpha) \;\mathcal{U}_m(n; \alpha, \beta, N) -n\;\mathcal{U}_m(n-1; \alpha, \beta, N)=-\alpha \; \mathcal{U}_m(n; \alpha +1, \beta, N),
\end{equation}
using the fact that $U_m(n)=\epsilon_m \;\mathcal{U}_m(n; \alpha, \beta, N)$ when the coefficients $\gamma^*_n$ given in  \eqref{gamma} are substituted in \eqref{Umn}.
With the explicit expression \eqref{Umn} of $\mathcal{U}_m(n; \alpha, \beta, N)$ as a ${_3}F_2(1)$, we may use the contiguity relation $X\;\mathcal{U}_m(n; \alpha, \beta, N)=-\alpha\; \mathcal{U}_m(n; \alpha, \beta, N)$ and its companion $Z\;\mathcal{U}_m(n; \alpha, \beta, N)=\frac{\alpha}{n-\alpha}\;\mathcal{U}_m(n; \alpha, \beta, N)$ that follows from \eqref{XZdiff}, to obtain the tridiagonal representation given in \eqref{X+}, \eqref{Xd}, \eqref{X-}, \eqref{Z+}, \eqref{ZD}, \eqref{Z-} of the algebra whose generators $Z$ and $X$ satisfy $[Z, X]=Z^2+Z$. It is also seen that
\begin{align}
    YU_m(n;\alpha)=& (e_m(x;\beta), V^TX^Td^*(x;\alpha))=(Ve_m(x;\beta), X^Td^*(x; \alpha) \nonumber \\
    =&-\alpha \nu_m (e_m(x;\beta), d^*(x;\alpha +1)=\alpha n(n-N+\beta)U_m(n; \alpha+1),
\end{align}
which using the expression \eqref{Ydiff} for $Y$ as a difference operator yields another contiguity relation for $\mathcal{U}_m(n; \alpha, \beta, N)$ (see \cite{tsujimoto2020algebraic}).
\subsubsection{The biorthogonal partners and orthogonality relations}
The biorthogonal partners to the functions $U_m(n)$ are the scalar products $\Tilde{U}_m(n)=(e^*_m, Zd_n)$. Using again the differential model we have
\begin{equation}
    \Tilde{U}_m(n)=\frac{1}{2\pi i} \oint_{\Gamma} dx\; \gamma_n \delta^*_m\; x^{-1+n-N} (1-x)^{1-\alpha}{_2}F_1 \left( {m-N, -m-\beta \atop  -N } ; x\right).
\end{equation}
Evaluating this integral amounts to identifying the coefficient of $x^{N-n}$ in the expansion
\begin{equation}
  (1-x)^{1-\alpha}{_2}F_1 \left( {m-N, -m-\beta \atop  -N } ; x\right)= \frac{1}{\delta^*_m}\sum_{\ell} \frac{1}{\gamma_{\ell}}
  \Tilde{U}_m (\ell) x^{N-\ell}\label{expand2}
\end{equation}
where again the domain of $\Tilde{U}_m $ has been formally extended in the right-hand side for the sake of notation. To achieve this task, we may rely on the following fact which is derived in Appendix \ref{B}. In view of the lower parameter $-N$, the polynomial ${_2}F_1 \left( {m-N, -m-\beta \atop  -N } ; x\right)$ is only defined for powers of $x$ up to $N$. This however suffices to identify the coefficient of $x^{N-n}$ in the left hand side of \eqref{expand2}. It so happens that although the functions $(1-x)^{\beta}{_2}F_1 \left( {-m, m-N+\beta \atop  -N } ; x\right)$ and ${_2}F_1 \left( {m-N, -m-\beta \atop  -N } ; x\right)$ are manifestly not equal (one being a polynomial the other not), they coincide term by term up to $x
^N$. Since this is the only range that matters for our computation we might hence replace one by the other and look alternatively for the coefficient of $x^{N-n}$ in
\begin{equation}
  (1-x)^{1-\alpha+\beta}{_2}F_1 \left( {-m, m-N+\beta \atop  -N } ; x\right)= \frac{1}{\delta^*_m}\sum_{\ell} \frac{1}{\gamma_{\ell}}
  \Tilde{U}_m (\ell) x^{N-\ell}. \label{expand3}
\end{equation}
This turn out to be simpler since it has already been performed. Indeed comparing with \eqref{expand1}, we see that $\frac{1}{\delta^*_m\gamma_{n}}
  \Tilde{U}_m (n)$ will be obtained from the expression for $\frac{1}{\delta_m\gamma^*_{n}} U_m (n)$ given in \eqref{Umn} by making the substitutions
  $n \rightarrow (N-n), \;\alpha \rightarrow (\beta+2-\alpha)$ (see \cite{tsujimoto2020algebraic}). We thus find
  \begin{equation}
 \Tilde{U}_m(n) = \delta^*_m \gamma_n
 \frac{(\alpha-\beta -1)_{N-n}}{(N-n)!}
 \; {_3}F_2 \left({ N-n,-m, m-N+\beta \atop -N, \beta +2 -\alpha -n }  ; 1 \right). \label{Umntilde}
\end{equation}
 Up to a normalization factor, this gives  the biorthogonal partner \cite{tsujimoto2017tridiagonal}
 \begin{equation}
   \mathcal{V}_m(n;\alpha,\beta,N)=\mathcal{U}_m(N-n;\beta+2-\alpha,\beta,N)   \label{dratHahn}
 \end{equation}
 of the rational Hahn function $\mathcal{U}_m(n;\alpha,\beta,N)$  multiplied by $\delta^*_m\gamma_n \frac{(\alpha-\beta -1)_{N-n}}{(N-n)!} $.

 We can now identify the weight function with respect to which these rational functions of Hahn type are orthogonal to one another. We have already established in Subsection \ref{sub} the biorthogonality of the functions $U_m(n)$ and $\Tilde{U}_k(n)$. Recall \eqref{wgg}. Omitting the factor $\delta_m \delta^*_k$ and other normalization constants, the biorthogonality relation $\sum_{n=0}^{N}\Tilde{U}_k(n)U_m(n)w_n=0 \; \text{if} \; k\neq m$ translates into
 \begin{equation}
     \sum_{n=0}^N \frac{(\alpha-\beta -1)_{N-n}}{(N-n)!} \frac{(1-\alpha)_n}{n!}\mathcal{V}_k(n;\alpha,\beta,N)\mathcal{U}_m(n;\alpha,\beta,N)=0\quad \text{if} \quad k\neq m \label{orel}
 \end{equation}
 using \eqref{Umn}, \eqref{Umntilde}, \eqref{ratHahn}, \eqref{dratHahn} and \eqref{wgg}. Note that the factor $\gamma_n \gamma^*_n$ is cancelled through the product with $w_n$. Using the identities \eqref{pochrel}, the first terms in \eqref{orel} can be transformed to make this equation read
 \begin{equation}
    \sum_{n=0}^N \mathcal{V}_k(n;\alpha,\beta,N)\mathcal{U}_m(n;\alpha,\beta,N) \mathfrak{w}^{(\alpha, \beta)}_n =0\quad \text{if} \quad k\neq m,
 \end{equation}
 with
 \begin{equation}
    \mathfrak{w}^{(\alpha, \beta)}_n = \frac{(\beta-\alpha-N+2)_N}{(\beta-N+1)_N} \frac{(-N)_n(1-\alpha)_n}{n!(\beta-\alpha-N+2)_n}.\label{weight}
 \end{equation}
 We thus recover the weight function given in \cite{tsujimoto2020algebraic} which intervenes in the orthogonality of the rational functions of Hahn type. The $n$-independent factors in \eqref{weight} have been introduced so that $\sum_{n=0}^N \mathfrak{w}^{(\alpha, \beta)}_n=1$.

 \subsection {The Hahn orthogonal polynomials}
 A central point we are making in this study is that the meta-Hahn algebra offers a unified description of both the rational Hahn functions and the Hahn polynomials. In this section we are using the differential model to confirm this and obtain explicitly the special functions. Having dealt with the rational functions, we now come to the orthogonal polynomials. These have been identified in Subsection \ref{5.2} with the scalar products $S_m(n)=(e_m, f^*_n)$
 and $\Tilde{S}_m(n)=(e^*_m, f_n)$ where $f_n$ and $f^*_n$ are respectively the eigenvectors of the linear pencil $W=X+\mu Z$ and its adjoint $W^T$ appropriately defined so that their common spectrum is (see Remark \ref{4.1})
\begin{equation}
    \rho_n=-n-\alpha-\mu.
\end{equation}
Using the differential realization for $X$ and $Z$ given in \eqref{X_D} and \eqref{Z_D} and the expressions \eqref{X_DT} and \eqref{Z_DT} for their Lagrange adjoints, we find the following solutions for the EVPs $Wf_n(x;\mu)= \rho_nf_n(x;\mu)$ and $W^Tf^*_n(x; \mu)=\rho_n f^*_n(x; \mu)$:
\begin{align}
    f_n(x;\mu)&=\epsilon_n (x-1)^{\mu}\left(\frac{x}{x-1}\right)^n,\\ \label{fn}
    f^*_n(x;\mu)&=-\epsilon^*_n(x-1)^{-\mu -2}\left(\frac{x}{x-1}\right)^{-n-1}, \qquad n=0, \dots, N,
\end{align}
where $\epsilon_n$ and $\epsilon^*_n$ are normalization constants. One readily verifies that the functions $f^*_k(x;\mu)$ and $f_n(x;\mu)$ are mutually orthogonal:
\begin{align}
    (f^*_k(x;\mu), f_n(x;\mu))=&-\epsilon^*_k\epsilon_n\frac{1}{2\pi i} \oint_{\Gamma} (1-x)^{-2} \left(\frac{x}{x-1}\right)^{n-k-1}dx \nonumber\\
    =&\;\epsilon^*_k\epsilon_n\frac{1}{2\pi i} \oint_{\Gamma '} y^{n-k-1} dy = \epsilon^*_n\epsilon_n \delta_{k,n}
\end{align}
using the change of variable
\begin{equation}
    y=\frac{x}{x-1}.\label{y}
\end{equation}
With $e_m(x;\beta)$ and $e^*_m(x;\beta)$ given by \eqref{ereal} and \eqref{e*real}, we may now proceed to determine $S_m(n)$ and $\Tilde{S}_m(n)$ as we have done for $U_m(n)$ and $\Tilde{U}_m(n)$. Consider first
\begin{align}
    S_m(n)&=(f^*_n(x;\mu), e_m(x;\beta)) \label{integ}\\ \nonumber
    &=-\epsilon^*_n\delta_m\frac{1}{2\pi i} \oint_{\Gamma}dx (x-1)^{-\mu -2}\left(\frac{x}{x-1}\right)^{-n-1}\; {_2}F_1 \left( {-m, m-N+\beta \atop  -N } ; x\right).
\end{align}
In order to perform this integral it is useful to call upon the identity
\begin{eqnarray}
{_2}F_1 \left( {-m, m-N+\beta \atop  -N } ; x\right) = (-1)^m(x-1)^m \: {_2}F_1 \left( {-m, -m-\beta \atop  -N } ; \frac{x}{x-1}\right) \label{hyp_tr}
\end{eqnarray}
and to use the change of variable \eqref{y} to transform the integral in \eqref{integ} into
\begin{equation}
    S_m(n)=(-1)^m\epsilon^*_n\delta_m\frac{1}{2\pi i} \oint_{\Gamma}dy (y-1)^{\mu -m} y^{-n-1}\; {_2}F_1 \left( {-m, -m-\beta \atop  -N } ; y\right).
\end{equation}
In a familiar fashion, we see that the function $S_m(n)$ will be provided by the coefficient of $y^n$ in the expansion
\begin{equation}
    (y-1)^{\mu-m}{_2}F_1 \left( {-m, -m-\beta \atop  -N } ; y\right)=\frac{(-1)^m}{\delta_m}\sum_{\ell}\frac{1}{\epsilon^*_{\ell}}S_m(\ell)y^{\ell} \label{expand4}
\end{equation}
where the domain of $S_m$ is extended on the right hand side to avoid additional notation. Using the binomial theorem \eqref{binom} and the explicit formula for the Gauss hypergeometric function, we find
\begin{eqnarray}
S_m(n)= (-1)^m \delta_m \epsilon^*_n \frac{(m-\mu)_n}{n!} \:  {_3}F_2 \left( {-m, -m-\beta, -n \atop  -N, \mu + 1 -m-n } ; 1 \right). \label{Smn}
\end{eqnarray}
At this point we bring the following transformation formula \footnote{This formula can be obtained from an example in \cite{bailey1964generalized} or from the Whipple formula \cite{gasper2004basic} when one upper and one lower parameter in the terminating ${_4}F_3(1)$ are set equal and sent to infinity.} for the terminating ${_3}F_2(1)$:
\begin{equation}
{_3}F_2 \left( {-n, b, c \atop  d, e} ; 1 \right) = \frac{(e-c)_n}{(e)_n} \: {_3}F_2 \left( {-n, d-b, c \atop  d, c-e-n+1} ; 1 \right). \label{tr32}
\end{equation}
With the identification
\begin{equation}
    b=m+\beta-N, \quad c=-m, \quad d=-N, \quad e=-\mu,
\end{equation}
we have
\begin{equation}
     S_m(n)=(-1)^m\delta_m\epsilon^*_n \frac{(-\mu)_n}{n!}\: {_3}F_2 \left( {-m, m+\beta-N, -n \atop  -N, -\mu} ; 1 \right) \label{662}
\end{equation}
and redefining the parameters
\begin{eqnarray}
\mu= -\hat{ \alpha} -1, \; \beta = \hat{\alpha}+\hat{\beta}+N+1  \label{tr_parH}
\end{eqnarray}
as in \eqref{ii} and \eqref{iii}, we arrive at
\begin{equation}
    S_m(n)=(-1)^m \delta_m \epsilon^*_n \frac{(\hat{ \alpha} +1)_n}{n!} Q_m(n, \hat{\alpha}, \hat{\beta}, N) \label{untils}
\end{equation}
where $Q_m(n, \hat{\alpha}, \hat{\beta}, N)$ are Hahn polynomials as per \eqref{HahnOP}.

In order to complete the picture, we shall examine the other related function: $\Tilde{S}_m(n)=(e^*_m, f_n)$. From \eqref{e*real} and \eqref{fn} we have
\begin{align}
    \Tilde{S}_m(n)&=(f_n(x;\mu), e^*_m(x;\beta)) \label{integ2}\\ \nonumber
    &=\epsilon_n\delta^*_m\frac{1}{2\pi i} \oint_{\Gamma}dx \;(x-1)^{\mu -1-N}\left(\frac{x}{x-1}\right)^{-1+n-N}\; {_2}F_1 \left( {m-N, -m-\beta \atop  -N } ; x\right).
\end{align}
Performing the change of variable $y=\frac{x}{x-1}$ one has
\begin{equation}
     \Tilde{S}_m(n)=
    \epsilon_n\delta^*_m\frac{1}{2\pi i} \oint_{\Gamma} dy \;(y-1)^{N-\mu -1}\; y^{-1+n-N}\; {_2}F_1 \left( {m-N, -m-\beta \atop  -N } ; \frac{y}{y-1}\right).
\end{equation}
To compute this integral we need the coefficient of $y^{N-n}$ in the expansion of the two other factors of the integrand. This implies that the monomials in $y$ of degree higher than $N$ in the hypergeometric series will not affect the result. We may then use the fact, also explained in Appendix B, that ${_2}F_1 \left( {m-N, -m-\beta \atop  -N } ; \frac{y}{y-1}\right)$ coincides with $(1-y)^{-m-\beta}\;{_2}F_1 \left( {-m, -m-\beta \atop  -N } ; y\right)$ up to the power $y
^N$ to replace one function by the other and look for the coefficient of $y^{N-n}$ in
\begin{equation}
    (y-1)^{N-\mu-m-\beta-1}{_2}F_1 \left( {-m, -m-\beta \atop  -N } ; y\right)=\frac{(-1)^{m+\beta+1}}{\delta^*_m}\sum_{\ell}\frac{1}{\epsilon_{N-\ell}}S_m(N-\ell)y^{\ell}.
\end{equation}
As it happens, we have carried out such an evaluation in \eqref{expand4}. Upon comparing we see that we only need to do the substitutions:
\begin{equation}
    \mu \;\rightarrow \; N-\mu-\beta -1 \qquad \text{and} \qquad n\;\rightarrow\; N-n
\end{equation}
in the expression for $\frac{(-1)^m }{\delta_m\epsilon^*_n}S_m(n)$ given in \eqref{662} to obtain the desired coefficient. This gives
\begin{equation}
     \Tilde{S}_m(n)=(-1)^{m+\beta+1}\delta^*_m\epsilon_n \frac{(\mu+\beta+1-N)_{N-n}}{(N-n)!}\: {_3}F_2 \left( {-m, m+\beta-N, n-N \atop  -N, \mu+\beta+1-N} ; 1 \right).
\end{equation}
Using the relations \eqref{tr_parH} we arrive at
\begin{equation}
     \Tilde{S}_m(n)=(-1)^{m+\beta+1}\delta^*_m\epsilon_n \frac{(\hat{\beta}+1)_{N-n}}{(N-n)!}\: {_3}F_2 \left( {-m, m+\hat{\alpha}+\hat{\beta}+1, n-N \atop  -N, \hat{\beta}+1} ; 1 \right).
\end{equation}
It can be seen from the difference equation obeyed by the Hahn polynomials that they have the following property under the exchange $n \rightarrow N-n$:
\begin{equation}
    Q_m(n,\hat{\alpha}, \hat{\beta}, N) = (-1)^m Q_m(N-n,\hat{\beta}, \hat{\alpha}, N).
\end{equation}
In light of this, we see that
\begin{equation}
   \Tilde{S}_m(n)=\hat{\delta}^*_m\epsilon_n \frac{(\hat{\beta}+1)_{N-n}}{(N-n)!}  Q_m(n,\hat{\alpha}, \hat{\beta}, N), \label{tilS}
\end{equation}
where signs have been absorbed in $\hat{\delta}^*_m$. We thus verify what we have previously established on general grounds: that $S_m(n)$ and $\Tilde{S}_m(n)$ involve the same Hahn polynomials.

As a by-product, we can readily obtain the weight function for the Hahn polynomials. We know from Subsection \ref{5.2} that $\sum_{n=0}^N \Tilde{S}_m(n)S_k(n)\zeta_n=0$ if $m \neq k$, where $(f^*_m,f _n)=\zeta ^{-1}_n \delta _{m,n}$. Here $\zeta ^{-1}_n=\epsilon^*_n\epsilon_n$. Using \eqref{untils} and \eqref{tilS} and dropping normalization factors, we have
\begin{equation}
    \sum_{n=0}^N\; \mathfrak{W}^{(\hat{\alpha}, \hat{\beta})}_n \; Q_m(n,\hat{\alpha}, \hat{\beta}, N) Q_k(n,\hat{\alpha}, \hat{\beta}, N) = 0 \quad \text{if} \quad m\neq k,
\end{equation}
where
\begin{equation}
  \mathfrak{W}^{(\hat{\alpha}, \hat{\beta})}_n = \frac{(\hat{ \alpha} +1)_n}{n!}  \frac{(\hat{\beta}+1)_{N-n}}{(N-n)!}.
\end{equation}
This expression coincides \cite{koekoek2010hypergeometric} (up to $n$-independent factors) with the standard weight function associated to the Hahn polynomials.

Note that the differential model could be further exploited to obtain additional properties of the special functions connected to the meta-Hahn algebra. It will not have escaped the reader's attention that this model has led to interesting representations of the Hahn rational functions and Hahn polynomials in terms of contour integrals.

 \section{Embedding of the meta-Hahn $m\mathfrak{H}$ algebra into the universal envelopping algebra  $\mathcal{U}(\mathfrak{sl}_2)$}
We shall here make the observation that $m\mathfrak{H}$ admits an embedding into $\mathcal{U}(\mathfrak{sl}_2)$. The Lie algebra $\mathfrak{sl}_2$ has three generators $\{J_0, J_{\pm}\}$ that satisfy the commutation relations
\begin{eqnarray}
[J_0,J_{\pm}]=\pm J_{\pm} , \quad [J_+,J_-]=2 J_0. \label{sl_2}
\end{eqnarray}
Its Casimir element $C \in \mathcal{U}(\mathfrak{sl}_2)$ is given by
\begin{eqnarray}
C = J_0^2-J_0 +J_+ J_-  . \label{Q_expr}
\end{eqnarray}
In irreducible representations $C$ is a multiple of the identity written in the form $C=j(j+1)\mathcal{I}$ and in finite-dimensional cases $j=\frac{N}{2}$ where $N+1$ is the dimension.

We shall provide below the map of the meta-Hahn algebra into (the completion of) $\mathcal{U}(\mathfrak{sl}_2)$. To offer first the results of the computation in general, we shall consider the relations
\begin{align}
[Z,X]=&Z^2+Z, \label{mhgzx}\\
[X,V] =& \{V,Z\}+ \eta_1 X + V + \eta_2 Z + \eta_0 \mathcal{I},\label{mhgxv}\\
[V,Z] =& \eta_4 X + \eta_1 Z + \eta_3 \mathcal{I},\label{mhgvz}
\end{align}
obtained from \eqref{ZXA}, \eqref{XVA}, \eqref{VZA} after the constraints \eqref{constr_meta} imposed by the Jacobi identity have been implemented.

\begin{prop}
The following formulas provide an embedding of the algebra defined by \eqref{mhgzx}, \eqref{mhgxv}, \eqref{mhgvz} into $\mathcal{U}(\mathfrak{sl_2)}$:
\begin{eqnarray}
&&Z= J_+ - \mathcal{I} , \quad X= -J_0 J_+ + J_0 + \xi_1 J_+  + \xi_0 \mathcal{I}  , \nonumber \\
&&V=  \xi_2 J_0^2 + \xi_3 J_0 + \xi_4 J_-  + \xi_5 \mathcal{I}
\label{meta_sl} \end{eqnarray}
with
\begin{eqnarray}
&&\xi_0= \frac{\eta_1-\eta_4}{\eta_3}, \quad \xi_2 = -\frac{\eta_4}{2}, \quad \xi_3 = \eta_1 +\eta_4(\xi_1 -1/2), \nonumber \\
&& \xi_4 = -\frac{\eta_4}{2}, \quad  \xi_5 = \frac{1}{2} \left(\eta_4(\xi_1 - \xi_1^2) + \eta_1 (1-2 \xi_1) -\eta_2  \right)
\end{eqnarray}
and the parameter $\xi_1$ a root of the quadratic equation
\begin{equation}
-{\eta_{4}}^{2}{\xi_{{1}}}^{2}+\eta_4 \left( -2\,\eta_1+\eta_4
 \right) \xi_{{1}} -2\,\eta_1^{2}+2\,\eta_1 \eta_3 -2
\,\eta_0 \eta_4 + \eta_4\eta_1 + \eta_2 \eta_4 + \eta_4^{2} \,C  =0.
\end{equation}
\end{prop}
The Casimir operator \eqref{Q_expr} of the (non-standardized) meta-Hahn algebra takes the following expression under this embedding:
\begin{equation}
Q= \eta_4 \, C -2\,\eta_0 + \eta_2 -{\frac { \left( \eta_3 - \eta_1 \right) ^{2}}{\eta_4}}. \label{Q_sl2}
\end{equation}

We now revert to the standardized meta-Hahn algebra where
\begin{equation}
    \eta_4=2 \qquad \text{and} \qquad \eta_2=-\eta_1
\end{equation}
and moreover set $\eta_3= 2\alpha -\beta -1$ as in \eqref{sv1} and take the central charges $\eta_0$ and $\eta_1$ to be  $\eta_0 = (N-1-\beta)\alpha + \beta +1$, $\eta_1=N-1-\beta$, as in \eqref{sv} in the $(N+1)$-dimensional representation where $C=\frac{N}{2}(\frac{N}{2}+1)$.
It is then found that the coefficients $\xi_i, \;i=0,\dots, 5$, take the simple expressions:
\begin{eqnarray}
\xi_0=\frac{N}{2}-\alpha,  \quad \xi_1 = 1-\frac{N}{2} , \quad \xi_2 = -1,  \quad
\xi_3 =-\beta, \quad \xi_4 =-1, \quad \xi_5 = \frac{N}{2}\left(\frac{N}{2}-\beta\right). \label{par_mHahn}
\end{eqnarray}
and that the value of the Casimir element \eqref{Q} is
\begin{eqnarray}
Q=-1+2\,\alpha-\beta+2\,\alpha\,\beta-2\,{\alpha}^{2}.
\end{eqnarray}

Let us remark to conclude this Section that the differential operators (again keeping the same notation for the abstract generators and their realizations):
\begin{align}
    J_0=&\;x\frac{d}{dx} + \tau \nonumber \\
    J_+=&\;x, \nonumber \\
    J_-=&\;-\left(x\frac{d^2}{dx^2}+2\tau\frac{d}{dx}\right),\label{sl2real}
\end{align}
verify the commutation relations \eqref{sl_2} of $\mathfrak{sl}_2$ with $\tau$ a parameter. Take
\begin{equation}
    \tau=-\frac{N}{2}
\end{equation}
(in which case $C=\frac{N}{2}(\frac{N}{2}+1)$), it is readily observed that for this choice of $\tau$, the differential realization \eqref{Z_D}, \eqref{X_D}, \eqref{V_D} of $m\mathfrak{H}$ at the heart of Section \ref{6} is recovered when the differential operators \eqref{sl2real} are substituted in the formulas \eqref{meta_sl} with the coefficients given by \eqref{par_mHahn}.

\section{Conclusion}

We want to stress in closing that the analysis we have presented sets the stage for a novel algebraic description of the polynomials of the Askey scheme coupled to associated biorthogonal rational functions. Before we comment on that let us first summarize the findings we reported.

We introduced an algebra $m\mathfrak{H}$ that we called the meta-Hahn algebra. The reason for the name is that it admits embeddings of both the Hahn algebra and the rational Hahn algebra which account for the bispectral properties of the Hahn orthogonal polynomials and rational functions respectively. The algebra $m\mathfrak{H}$ thus unifies the algebraic description of these two sets of functions. The algebra $m\mathfrak{H}$ has three generators $Z, X, V$ subjected to relations that involve at most one quadratic term. Its construction is centered around the subalgebra generated by $Z$ and $X$ which is a real form of the deformed Jordan plane. With an eye to the connection to special functions, certain elements of the representation theory of $m\mathfrak{H}$ were developed. To that end, various bases of a finite-dimensional module were considered, namely the one associated to the GEVP defined by $Z$ and $X$, the eigenbases of $V$ and of the linear pencil $X+\mu Z$ and the corresponding adjoint bases. It was shown in particular that $Z$ and $X$ act tridiagonally in the eigenbasis of $V$ and that $V$ is a Hessenberg matrix in the GEVP basis. The features identified have allowed to establish that the overlaps between these bases are biorthogonal rational functions or orthogonal polynomials that are bispectral. (It would be interest in the future to complete the study of the representation theory of $m\mathfrak{H}$.) A realization of $m\mathfrak{H}$ in terms of differential operators was presented and used to obtain the explicit expressions for the functions defined by these overlaps. The Hahn polynomials and rational functions were thus found and the model provided much of their characterization ((bi-)orthogonality, contiguity relations, etc.) leading in particular to representations of these functions in terms of contour integrals. It should be underscored moreover that the subalgebra generated solely by $Z$ and $X$ is on its own intimately connected to the Hahn polynomials. Indeed its representations by tridiagonal matrices (which we have in the eigenbasis of $V$) were seen to entail the recurrence relations of these orthogonal polynomials. It was noted additionally that $m\mathfrak{H}$ derives from a potential and admits an embedding in  $\mathcal{U}(\mathfrak{sl}_2)$. While the appendices that follow contain technical details, we want to point out that the computations in Section \ref{6} led to the identification of a Pad\'e approximant for the binomial series in terms of Jacobi polynomials with a negative integer. This result which will be found in Appendix \ref{AppB} stems from traces of the Euler and Pfaff transformations of the hypergeometric series which per se are not valid when the denominator parameter is a negative integer.

Let us now explain why we believe that this study is paving the way for a new algebraic picture of the Askey scheme extended in fact to include classes of biorthogonal rational functions. The meta-Hahn algebra can be viewed as an enlargement through the addition of the generator $V$ of the non-commutative model of the plane given in this case by the algebra generated by $Z$ and $X$ subjected to $[Z,X]=Z^2+Z$. We have noted the connection of this algebra with the Hahn polynomials. The commutation relations of $V$ with $X$ and $Z$ are minimally quadratic and such that the resulting algebra admits an embedding of the Hahn algebra. Now, there are a number of two-generated algebras and a classification can be found in \cite{gaddis2015two}, see also \cite{smith1992quantum}. Such algebras have arisen also in studies of martingale polynomials \cite{bryc2007quadratic} and in the physical context of asymmetric exclusion models \cite{essler1996representations}. We trust that these are related to hypergeometric polynomials like our real version of the deformed Jordan plane. It is hence our conviction that algebras analogous to the meta-Hahn algebras can be defined as extensions of these other two-dimensional algebras and that they will henceforth provide descriptions of families of orthogonal polynomials and rational functions in a unified way. Let us mention the following in support of this. The $q$-oscillator algebra with generators $X$ and $Z$ verifying $XZ-qZX=1$ is obviously one of the relevant two-generated algebras. Its tridiagonal representations have been obtained in \cite{tsujimoto2017tridiagonal} and found to yield the recurrence coefficients of the Askey--Wilson polynomials. The analog of the generator $V$ was also found and the resulting algebra was identified as that of the big $q$-Jacobi polynomials. The exploration of the associated rational function was not carried out and this is obviously on our to-do list now. Another two-generated algebra of interest is the one with relation $[Z,X]=Z
^2+X$ which is discussed in \cite{gaddis2015two}. Preliminary investigations are indicating that this algebra is related to the Racah polynomials. We have therefore good reasons to believe that the construction of an algebra analogous to the meta-Hahn algebra and centered around this ``Racah" algebra will yield a framework offering an integrated picture for the Racah polynomials and rational functions of the ${_4}F_3$-type. We are poised to pursue this program and hope to report on it in the near future.

\appendix

\section{Orthogonality of the functions $e_n(x;\beta)$ and $e^*_m(x;\beta)$} \label{App}

We give in this Appendix the details of the computations showing that the functions $e_n(x;\beta)$ and $e^*_m(x;\beta)$ defined in \eqref{ereal} and \eqref{e*real} are mutually orthogonal.
The scalar product on $\mathfrak{V}_N$ has been defined by integration over the contour $\Gamma$ taken to be a circle $|x|= a<1$ in the complex plane, i.e.
\begin{equation}
    (f(x), g(x))=\frac{1}{2\pi i} \oint_{\Gamma} f(x)g(x)dx, \quad f, g \in \mathfrak{V}_N.
\end{equation}
We have
\begin{align}
(e^*_m(x;\beta),& e_n(x;\beta))= \nonumber \\
& \delta^*_m\delta_n\frac{1}{2\pi i} \oint_{\Gamma}dx\;x^{-1-N}{_2}F_1 \left( {m-N, -m-\beta \atop  -N } ; x\right){_2}F_1 \left( {-n, n-N+\beta \atop  -N } ; x\right)
\end{align}
and want to show that $(e^*_m(x;\beta), e_n(x;\beta))=0$ if $m\; \ne\; n $.
Clearly, this property is equivalent to the statement that in the power expansion
\begin{eqnarray}
\Psi_{nm}(x)= {_2}F_1\left( {-m, m+\beta-N  \atop -N};x \right) \ {_2}F_1\left( {n-N, -n-\beta  \atop -N};x \right) = \sum_{s=0}^{m-n+N} A_s(n,m) x^s
\end{eqnarray}
the coefficient of $x^N$ vanishes when $n \ne m$. We must therefore prove that
\begin{equation}
A_N(n,m) =0, \quad n \ne m.
\end{equation}
In order to show this, we consider product of the two terminating hypergeometric series in $x$:
\begin{equation}
\Psi_{nm}(x) = \sum_{s,i} \frac{(n-N)_s (-n-\beta)_s (-m)_i  (m+\beta-N)_i}{s! i! (-N)_s (-N)_i} x^{s+i}
\end{equation}
to find that
\begin{equation}
A_N(n,m)= \frac{(n-N)_{N-i} (-n-\beta)_{N-i} (-m)_i  (m+\beta-N)_i}{(N-i)! i! (-N)_{N-i} (-N)_i}.
\end{equation}
It is easily seen that the following product $(n-N)_{N-i} (-m)_i$ of Pochhammer symbols is zero for all $i=0,1,\dots, N$ if $n \ne m$. When $n=m$ there is one and exactly one value of $i$ for which this product is nonzero, it is $i=n$ and then
\begin{equation}
(n-N)_{N-i} (-n)_i = -n!(N-n)! \: \delta_{ni}
\end{equation}
We have thus have proven that $A_N(n,m)=0$ if $n \ne m$ and that $A_N(n,n) \ne 0$.

\section{Restricted versions of Euler's and Pfaff's transformations and a Pad\'e approximation table for the binomial function}\label{AppB}

This Appendix confirms the validity of the substitutions that were performed in Subsubsection 6.2.3 and Subsection 6.3 when computing $\Tilde{U}_m(n)$ and $\Tilde{S}_m(n)$. It offers restricted versions of Euler's and Pfaff's transformations of the (truncating) hypergeometric series when the lower parameter is a negative integer. A Pad\'e approximant of the binomial function is moreover obtained as a by-product.

The Euler transformation of the hypergeometric series \label{B}
\begin{equation}
    {_2}F_1\left( {a, b  \atop c};x \right) = (1-x)^{c-a-b} \: {_2}F_1\left( {c-a, c-b  \atop c};x \right)
\end{equation}
is not valid if $c=-N$, with $N$ a nonnegative integer. Indeed, in this case the hypergeometric series has singular terms.
It is nevertheless possible to derive a restricted version of this transformation that has proved most useful in the context of the differential model for $m\mathfrak{H}$. We have:
\begin{prop}\label{propB1}
Assume that $N,n$ are positive integers such that $n<N$ and that $b$ is a real non-integer parameter.  Then the  relation
\begin{equation}
{_2}F_1\left( {-n, b  \atop -N};x \right) = (1-x)^{-N+n-b} {_2}F_1\left( {n-N, -N-b  \atop -N};x \right) \label{B2}
\end{equation}
holds up to terms of degree $N$ in $x$.
\end{prop}
The proof proceeds from multiplying the binomial series
\begin{equation}
 (1-x)^{\beta} =\sum_{s=0}^{\infty} \frac{(-\beta)_s}{s!} x^s, \quad \beta=-N+n-b
\end{equation}
and the hypergeometric series
\begin{equation}
  {_2}F_1\left( {n-N, -N-b  \atop -N};x \right) = \sum_{s=0}^{N-n} \frac{(n-N)_s (-N-b)_s}{s! (-N)_s} x^s.
\end{equation}
Writing
\begin{equation}
(1-x)^{-N+n-b} {_2}F_1\left( {n-N, -N-b  \atop -N};x \right) = \sum_{k=0}^{\infty} A_k x^k,
\end{equation}
the coefficients $A_k$ are found to be
\begin{equation}
    A_k=\frac{(N-n+b)_k}{k!}\;{_3}F_2\left( {-k, n-N, -N-b  \atop -N, 1-k-N+n-b};x \right).\label{B6}
\end{equation}
Clearly these $A_k$ are only defined for $k \leq N$. The
balanced hypergeometric  ${_3}F_2(1)$ series occurring in \eqref{B6} can then be evaluated using the Saalsch\"utz formula
\begin{equation}
  {_3}F_2\left( {-k, \;\alpha, \;\beta  \atop \gamma, \;1+\alpha+\beta+\gamma-k};\;x \right) = \frac{(\gamma-\alpha)_k(\gamma-\beta)_k}{(\gamma)_k (\gamma-\alpha-\beta)_k}  \label{saal}
\end{equation}
and one thus establishes the validity of the proposition. Note that ${_2}F_1\left( {-n, b  \atop -N};x \right)$
on the left hand side of \eqref{B2} is a polynomial in $x$ of degree $n$ while on the right hand side $(1-x)^{-N+n-b} {_2}F_1\left( {n-N, -N-b  \atop -N};x \right)$ is the product of a polynomial of degree $N-n$ with the binomial (transcendental) function $(1-x)^{-N+n-b}$ and is therefore an infinite series in $x$. Nevertheless, truncating this series after the first $N+1$ gives the polynomial $A_0+A_1x+A_2x^2+ \dots,+ A_Nx^N$ that coincides with the one on the left hand side. This is the precise meaning of Proposition \ref{propB1}.

A similar result that was used to compute the functions $\Tilde{S}_m(n)$ arises when considering Pfaff's transformation in the same context. We can formulate:
\begin{prop}
Assume that $N,n$ are positive integers such that $n<N$ and that $b$ is a real non-integer parameter.  The  relation
\begin{equation}
(1-x)^b \: {_2}F_1\left( {-n, b  \atop -N};x \right) =  {_2}F_1\left( {n-N, b  \atop -N}; \frac{x}{x-1} \right)    \label{2nd}
\end{equation}
applies term by term as an identity between power series of $x$ until the monomial $x^N$ is reached.
\end{prop}
This proposition can be considered as a restricted version of Pfaff's transformation formula:
\begin{equation}
 (1-x)^b \: {_2}F_1\left( {a, b  \atop c};x \right) =  {_2}F_1\left( {c-a, b  \atop c}; \frac{x}{x-1} \right)
\end{equation}
which holds identically if $c \ne -N$.
The proof is parallel to that of Proposition \ref{propB1}. One expands the left and right hand sides of \eqref{2nd} in power of $x$ and uses the  Saalsch\"utz formula \eqref{saal} to observe the equality of both sides as long as the (resummed) coefficients are defined.
There is a striking application of Proposition \ref{propB1} to the Pad\'e approximation table of the binomial function. It can be formulated as follows:
\begin{prop}
Let $\beta$ be an arbitrary non-integer parameter. The first $N+1$ terms of the binomial series
\begin{equation}
 (1-x)^{\beta} =\sum_{s=0}^{\infty} \frac{(-\beta)_s}{s!} x^s
\end{equation}
coincide with the first $N+1$ terms of the power expansion of the rational function
\begin{equation}
 R(x) = \frac{{_2}F_1\left( {n-N, -n-\beta  \atop -N};x \right)}{{_2}F_1\left( {-n, n+\beta-N  \atop -N};x \right)} = \sum_{s=0}^{N} \frac{(-\beta)_s}{s!} x^s + O(x^{N+1}).
\end{equation}
\end{prop}
This proposition is an immediate corollary of Proposition \ref{propB1}.
It allows to construct the whole Pad\'e interpolation table (see for example \cite{baker1996pade}) for the binomial function. Namely, we have that for every non-negative integers $n,m$, the binomial series coincides with the first $n+m+1$ terms of the power expansion of the rational function
\begin{equation}
  R_{mn}(x) = \frac{{_2}F_1\left( {-m, -n-\beta  \atop -n-m};x \right)}{{_2}F_1\left( {-n, -m+\beta  \atop -n-m};x \right)}.
\end{equation}
Hence the set of the rational functions $R_{mn}(x)$ for all $m,n=0,1,2,\dots$ constitutes the Pad\'e interpolation table for the binomial function $(1-x)^{\beta}$.


\section*{Acknowledgments}
The authors express their thanks to Andr\'e Beaudoin, Geoffroy Bergeron, Antoine Brillant, Maxim Derevyagin, Julien Gaboriaud and Satoshi Tsujimoto for enlightening discussions. The work of LV is supported in part by a Discovery Grant from the Natural Sciences and Engineering Research Council (NSERC) of Canada. AZ is gratefully holding Simons CRM Professorship and is funded by the National Foundation of China (Grant No.11771015).


\end{document}